\documentclass[a4paper,12pt,twoside]{article}
\usepackage[francais]{babel}
\usepackage{amssymb, amsmath, eucal}
\usepackage[all]{xy}
\usepackage{mathrsfs}
\usepackage[dvips]{graphics}
\begin{document}
\textwidth15.5cm
\textheight22.5cm
\voffset=-13mm
\newtheorem{The}{Theorem}[section]
\newtheorem{Lem}[The]{Lemma}
\newtheorem{Prop}[The]{Proposition}
\newtheorem{Cor}[The]{Corollary}
\newtheorem{Rem}[The]{Remark}
\newtheorem{Obs}[The]{Observation}
\newtheorem{SConj}[The]{Standard Conjecture}
\newtheorem{Titre}[The]{\!\!\!\! }
\newtheorem{Conj}[The]{Conjecture}
\newtheorem{Question}[The]{Question}
\newtheorem{Prob}[The]{Problem}
\newtheorem{Def}[The]{Definition}
\newtheorem{Not}[The]{Notation}
\newtheorem{Claim}[The]{Claim}
\newtheorem{Conc}[The]{Conclusion}
\newtheorem{Ex}[The]{Example}
\newtheorem{Fact}[The]{Fact}
\newcommand{\C}{\mathbb{C}}
\newcommand{\R}{\mathbb{R}}
\newcommand{\N}{\mathbb{N}}
\newcommand{\Z}{\mathbb{Z}}
\newcommand{\Q}{\mathbb{Q}}
\newcommand{\Proj}{\mathbb{P}}

\begin{center}

{\Large\bf Aeppli Cohomology Classes Associated with Gauduchon Metrics on Compact Complex Manifolds}

\end{center}

\begin{center}

 {\large Dan Popovici}

\end{center}

\vspace{1ex}

\noindent {\small {\bf Abstract.} We propose the study of a Monge-Amp\`ere-type equation in bidegree $(n-1,\,n-1)$ rather than $(1,\,1)$ on a compact complex manifold $X$ of dimension $n$ for which we prove ellipticity and uniqueness of the solution subject to positivity and normalisation restrictions. Existence will hopefully be dealt with in future work. The aim is to construct a special Gauduchon metric uniquely associated with any Aeppli cohomology class of bidegree $(n-1,\,n-1)$ lying in the Gauduchon cone of $X$ that we hereby introduce as a subset of the real Aeppli cohomology group of type $(n-1,\,n-1)$ and whose first properties we study. Two directions for applications of this new equation are envisaged\!: to moduli spaces of Calabi-Yau $\partial\bar\partial$-manifolds and to a further study of the deformation properties of the Gauduchon cone beyond those given in this paper.}

\section{Introduction}\label{section:introd}

 Let $X$ be a compact complex manifold, $\mbox{dim}_{\C}X=n$. The main theme of this paper is the interaction between various kinds of metrics (especially Gauduchon metrics) on $X$ and certain cohomology theories (especially the Aeppli cohomology) often considered on $X$. 

 On the metric side, let $\omega>0$ be a $C^{\infty}$ positive definite $(1,\,1)$-form (i.e. a Hermitian metric) on $X$. The following diagram sums up the definitions of well-known kinds of Hermitian metrics and the implications among them.

\vspace{3ex}

\noindent$\begin{array}{lllll}d\omega=0 & \Longrightarrow & \exists\,\, \alpha^{0,\,2}\in C^{\infty}_{0,\,2}(X,\,\C)\,\, \mbox{s.t.} & \Longrightarrow & \partial\bar\partial\omega=0 \\
 &  & d(\overline{\alpha^{0,\,2}}+\omega+\alpha^{0,\,2})=0 & & \\
 (\omega\,\,\mbox{is K\"ahler}) &   & (\omega\,\,\mbox{is Hermitian-symplectic}) &  & (\omega\,\,\mbox{is pluriclosed}) \\
\rotatebox{-90}{$\implies$} &  & & & \\
d\omega^{n-1}=0 & \Longrightarrow & \exists\,\, \Omega^{n-2,\,n}\in C^{\infty}_{n-2,\,n}(X,\,\C)\,\, \mbox{s.t.}  & \Longrightarrow & \partial\bar\partial\omega^{n-1}=0 \\
&  &  d(\overline{\Omega^{n-2,\,n}}+\omega^{n-1}+\Omega^{n-2,\,n})=0 & &   \\
(\omega\,\,\mbox{is balanced}) &   & (\omega\,\,\mbox{is strongly Gauduchon (sG)}) &  & (\omega\,\,\mbox{is Gauduchon}).\end{array}$

\vspace{2ex}

\noindent Recall that of the above six types of metrics, only Gauduchon metrics always exist ([Gau77]). Compact complex manifolds $X$ carrying any of the other five types of metrics inherit the name of the metric (K\"ahler, balanced, etc). 

 The conditions on the top line in the above diagram bear on the metric in bidegree $(1,\, 1)$, while those on the bottom line are their analogues in bidegree $(n-1,\, n-1)$. It is a well-known fact in linear algebra (see e.g. [Mic82]) that for every smooth $(n-1,\, n-1)$-form $\Gamma>0$ on $X$, there exists a unique smooth $(1,\, 1)$-form $\gamma>0$ on $X$ such that $\gamma^{n-1}=\Gamma$. We denote it $\gamma=\Gamma^{\frac{1}{n-1}}$ and call $\gamma$ the $(n-1)^{st}$ root of $\Gamma$. The power-root bijection between positive definite $C^{\infty}$ forms of types $(1,\,1)$ and $(n-1,\,n-1)$ suggests a possible duality between the metric properties in these two bidegrees. The following observation (noticed before in [IP02] and references therein as a consequence of more general results) gives a further hint. We give below a quick proof.

\begin{Prop}\label{Prop:pluri+bal} If $\omega$ is both pluriclosed and balanced, then $\omega$ is K\"ahler.

\end{Prop}

\noindent {\it Proof.} The pluriclosed assumption on $\omega$ translates to any of the following equivalent properties\!\!:

\begin{eqnarray}\label{eqn:pluriclosed-equiv}\partial\bar\partial\omega=0 \Longleftrightarrow \partial\omega\in\ker\bar\partial \Longleftrightarrow \star\,(\partial\omega)\in\ker\partial^{\star},\end{eqnarray}

\noindent the last equivalence following from the well known formula $\partial^{\star}=-\star\bar\partial\star$, where $\star=\star_{\omega}:\Lambda^{p,\,q}T^{\star}X\rightarrow \Lambda^{n-q,\,n-p}T^{\star}X$ is the Hodge-star isomorphism defined by $\omega$ for arbitrary $p,q=0,\dots , n$.

\noindent On the other hand, the balanced assumption on $\omega$ translates to any of the following equivalent properties\!\!:

\begin{eqnarray}\label{eqn:bal-equiv}d\omega^{n-1}=0 \Longleftrightarrow \partial\omega^{n-1}=0 \Longleftrightarrow \partial\omega \,\,\mbox{is primitive},\end{eqnarray}

\noindent the last equivalence following from the formula $\partial\omega^{n-1} = (n-1)\,\omega^{n-2}\wedge\partial\omega$. Now, since $\partial\omega$ is primitive by (\ref{eqn:bal-equiv}), a well-known formula (valid for any Hermitian metric $\omega$) given e.g. in [Voi02, Proposition 6.29, p. 150] spells

\vspace{1ex}

\hspace{10ex} $\star\,(\partial\omega) = i\,\frac{\omega^{n-3}}{(n-3)!}\wedge\partial\omega = \frac{i}{(n-2)!}\,\partial\omega^{n-2}.$ 

\vspace{1ex}

\noindent Since $\star\,(\partial\omega)\in\ker\partial^{\star}$ by (\ref{eqn:pluriclosed-equiv}), we get by applying $\partial^{\star}$ in the above identities\!\!:

\vspace{1ex}

\noindent $\partial^{\star}\partial\omega^{n-2}=0, \hspace{2ex}\mbox{so}\hspace{2ex} 0 = \langle\langle\partial^{\star}\partial\omega^{n-2},\, \omega^{n-2}\rangle\rangle = ||\partial\omega^{n-2}||^2, \hspace{2ex}\mbox{so}\hspace{2ex}\partial\omega^{n-2}=0.$

\vspace{1ex}

\noindent Since $\partial\omega^{n-2}=(n-2)\,\omega^{n-3}\wedge\partial\omega$, the last identity above gives $\omega^{n-3}\wedge\partial\omega=0$. Now, $\partial\omega$ is a form of degree $3$ while the Lefschetz operator on $3$-forms

$$L_{\omega}^{n-3} : \Lambda^3T^{\star}X\longrightarrow\Lambda^{2n-3}T^{\star}X, \hspace{3ex} \alpha\mapsto\omega^{n-3}\wedge\alpha,$$ 

\noindent is an isomorphism (see e.g. [Voi02, Lemma 6.20, p.146] -- no assumption on the Hermitian metric $\omega$ is needed). It follows that $\partial\omega=0$, i.e. $\omega$ is K\"ahler.   \hfill  $\Box$

\vspace{3ex}

 On the cohomological side, recall that for all $p, q=0, \dots , n$, the Bott-Chern cohomology group of type $(p,\,q)$ is defined as

$$H^{p, \, q}_{BC}(X, \, \C)=\frac{\ker(d:C^{\infty}_{p, \, q}(X)\rightarrow C^{\infty}_{p+1, \, q}(X) + C^{\infty}_{p,\,q+1}(X))}{\mbox{Im}(\partial\bar\partial:C^{\infty}_{p-1, \, q-1}(X)\rightarrow C^{\infty}_{p, \, q}(X))},$$

\noindent while the Aeppli cohomology group of type $(p,\,q)$ is defined as

$$H^{p, \, q}_A(X, \, \C)=\frac{\ker(\partial\bar\partial:C^{\infty}_{p, \, q}(X)\rightarrow C^{\infty}_{p+1, \, q+1}(X))}{\mbox{Im}(\partial:C^{\infty}_{p-1, \, q}(X)\rightarrow C^{\infty}_{p, \, q}(X)) + \mbox{Im}(\bar\partial:C^{\infty}_{p, \, q-1}(X)\rightarrow C^{\infty}_{p, \, q}(X))}.$$

\noindent There are always well-defined linear maps from $H^{p, \, q}_{BC}(X, \, \C)$, from $H^{p, \, q}_{\bar\partial}(X, \, \C)$ (the Dolbeault cohomology group of type $(p,\,q)$) and from $H^{p+q}(X,\,\C)$ (the De Rham cohomology group of degree $p+q$) to $H^{p, \, q}_A(X, \, \C)$ but, in general, they are neither injective, nor surjective. 

 We will be often considering the case when $X$ is a $\partial\bar\partial$-manifold. This means that the $\partial\bar\partial$-lemma holds on $X$, i.e. for all $p, q$ and for any smooth $d$-closed form $u$ of pure type $(p,\,q)$ on $X$, the conditions of $d$-exactness, $\partial$-exactness, $\bar\partial$-exactness and $\partial\bar\partial$-exactness are all equivalent for $u$.  

 If $X$ is a $\partial\bar\partial$-manifold, $H^{p, \, q}_A(X, \, \C)$ is canonically isomorphic to each of the vector spaces $H^{p, \, q}_{BC}(X, \, \C)$ and $H^{p, \, q}_{\bar\partial}(X, \, \C)$, while injecting canonically into $H^{p+q}(X,\,\C)$ (cf. Theorem \ref{The:ddbar-Aeppli}). In particular, if $(X_t)_{t\in\Delta}$ is a deformation of the complex structure of $X=X_0$, the various Aeppli cohomology groups of the fibres $X_t$ depend on $t$ but, if $X_0$ is assumed to be a $\partial\bar\partial$-manifold (in which case every $X_t$ with $t$ sufficiently close to $0$ is again a $\partial\bar\partial$-manifold by Wu's main theorem in [Wu06]), then for each $(p,\, q)$, all the groups $H^{p, \, q}_A(X_t, \, \C)$ inject canonically into a fixed De Rham cohomology group of $X$\!\!:

$$H^{p, \, q}_A(X_t, \, \C)\hookrightarrow H^{p+q}(X,\,\C), \hspace{3ex} t\in\Delta,$$

\noindent after possibly shrinking $\Delta$ about $0$. Under the same $\partial\bar\partial$ assumption on $X_0$ (hence also on $X_t$ for $t$ close to $0$), there are canonical isomorphisms (cf. Theorem \ref{The:ddbar-Aeppli})\!\!:

$$H^k(X,\,\C)\simeq\bigoplus\limits_{p+q=k}H^{p, \, q}_A(X_t, \, \C), \hspace{3ex} t\in\Delta,\hspace{1ex} k=0,\dots, 2n.$$

\noindent They depend only on the complex structure of $X_t$ and will be called the {\bf Hodge-Aeppli decomposition} of $X_t$ for $t$ in a possibly shrunk $\Delta$.

\vspace{3ex}

 We now bring together the metric and the cohomological points of view. Let $\omega$ be a Gauduchon metric on the $\partial\bar\partial$-manifold $X=X_0$, i.e. a $C^{\infty}$ positive definite $(1,\, 1)$-form such that $\partial\bar\partial\omega^{n-1}=0$. Then $\omega^{n-1}$ defines an Aeppli cohomology class $[\omega^{n-1}]_A\in H^{n-1, \, n-1}_A(X, \, \C)$ that we call the induced {\bf Aeppli-Gauduchon class}. The image $\{\omega^{n-1}\}\in H^{2n-2}(X,\,\C)$ under the canonical injection $H^{n-1, \, n-1}_A(X, \, \C)\hookrightarrow H^{2n-2}(X,\,\C)$ induced by the $\partial\bar\partial$ assumption on $X$ of the Aeppli-Gauduchon class $[\omega^{n-1}]_A$ will be called the associated {\bf De Rham-Gauduchon class}. Note that $\omega^{n-1}$ need not be $d$-closed, hence it need not define directly a De Rham class, but we have just argued that on a $\partial\bar\partial$-manifold $X$ there is a De Rham class of degree $2n-2$ (that we denote a bit abusively by $\{\omega^{n-1}\}\in H^{2n-2}(X,\,\C)$) canonically associated with the Aeppli class $[\omega^{n-1}]_A\in H^{n-1, \, n-1}_A(X, \, \C)$.

 Extending the approach of [Pop13a] from balanced classes to Gauduchon classes, we can define the fibres that are {\bf co-polarised} by the De Rham-Gauduchon class $\{\omega^{n-1}\}\in H^{2n-2}_{DR}(X,\,\C)$ in the family $(X_t)_{t\in\Delta}$ as being those $X_t$ for which $\{\omega^{n-1}\}$ remains of type $(n-1,\, n-1)$ in the Hodge-Aeppli decomposition 

\begin{equation}\label{eqn:Hodge-Aeppli-2n-2}H^{2n-2}(X,\,\C)\simeq H^{n, \, n-2}_A(X_t, \, \C)\oplus H^{n-1, \, n-1}_A(X_t, \, \C)\oplus H^{n-2, \, n}_A(X_t, \, \C)\end{equation}

\noindent of degree $2n-2$ on $X_t$, i.e. those $X_t$ for which the components of $X_t$-types $(n,\,n-2)$ and $(n-2,\,n)$ of $\{\omega^{n-1}\}\in H^{2n-2}(X,\,\C)$ vanish. (Since the class $\{\omega^{n-1}\}$ is real, it actually suffices for the $(n-2,\,n)$-component of $\{\omega^{n-1}\}$ to vanish.)

 We can construct a local deformation theory of Calabi-Yau $\partial\bar\partial$-manifolds co-polarised by a De Rham-Gauduchon class on the model of that for co-polarisations by a balanced class constructed in [Pop13a]. 

\vspace{2ex}

\noindent {\bf A Monge-Amp\`ere-type equation in bidegree $(n-1,\,n-1)$} 

\vspace{1ex}

 To go from local deformations to moduli spaces, we need canonical objects, namely we would like to single out in any co-polarising Gauduchon class $[\omega^{n-1}]_A$ (or $\{\omega^{n-1}\}$) a unique $(n-1)^{st}$ power of a Gauduchon metric for which we have prescribed the volume form. On a Calabi-Yau manifold $X$ (i.e. one for which the canonical bundle $K_X$ is trivial), this would entail the existence of a unique Ricci-flat Gauduchon metric $\omega$ of a certain shape whose Aeppli cohomology class $[\omega^{n-1}]_A\in H^{n-1,\,n-1}_A(X,\,\C)$ has been prescribed (arbitrarily). (By $\omega$ being Ricci-flat we mean that the Ricci form $\mbox{Ric}\,\omega$ of $\omega$ -- defined as the curvature form of the anti-canonical bundle $-K_X$ equipped with the metric induced by $\omega$ -- vanishes identically.)

 Motivated by considerations of this nature, we undergo to study in this and future work to which extent there is an Aeppli-Gauduchon analogue of Yau's theorem on the Calabi conjecture. Every representative of $[\omega^{n-1}]_A$ is of the form $\omega^{n-1} + \partial u + \bar\partial v$ with $u$ of type $(n-2,\, n-1)$ and $v$ of type $(n-1,\, n-2)$. To avoid an underdetermined equation, it seems sensible to look for forms of the special shape $u=\bar\partial\varphi\wedge\omega^{n-2}$ and $v=\partial\varphi\wedge\omega^{n-2}$ (up to constant factors), where $\varphi$ is a real smooth function on $X$ that we wish to find. We are thus led to look for positive definite $(n-1,\, n-1)$-forms that are Aeppli-cohomologous to $\omega^{n-1}$ of the shape \\

$\omega^{n-1} + \frac{i}{2}\,\partial(\bar\partial\varphi\wedge\omega^{n-2}) - \frac{i}{2}\,\bar\partial(\partial\varphi\wedge\omega^{n-2})$ \\

\hspace{20ex} $ = \omega^{n-1} + i\partial\bar\partial\varphi\wedge\omega^{n-2} + \frac{i}{2}\,(\partial\varphi\wedge\bar\partial\omega^{n-2} - \bar\partial\varphi\wedge\partial\omega^{n-2}).$

\vspace{2ex}

 Equations $(\star)$ and (\ref{eqn:main-eq-const}) proposed below involve taking the $(n-1)^{st}$ root of a positive definite $(n-1,\, n-1)$-form and thus produce a Gauduchon metric $\gamma$ with prescribed volume form $\gamma^n$ such that $\gamma^{n-1}$ is Aeppli-cohomologous to the $(n-1)^{st}$ power $\omega^{n-1}$ of the given Gauduchon metric $\omega$.

\begin{Question}\label{Question:equation} Let $X$ be a compact complex manifold of complex dimension $n\geq 2$. Fix an arbitrary Gauduchon metric $\omega$ on $X$. Consider the equation

\begin{equation}\nonumber\bigg[\bigg(\omega^{n-1} + i\partial\bar\partial\varphi\wedge\omega^{n-2} + \frac{i}{2}\,(\partial\varphi\wedge\bar\partial\omega^{n-2} - \bar\partial\varphi\wedge\partial\omega^{n-2})\bigg)^{\frac{1}{n-1}}\bigg]^n = e^f\,\omega^n\hspace{1ex}(\star)\end{equation}

\noindent subject to the positivity and normalisation conditions

\begin{equation}\label{eqn:conditions}\omega^{n-1} + i\partial\bar\partial\varphi\wedge\omega^{n-2} + \frac{i}{2}\,(\partial\varphi\wedge\bar\partial\omega^{n-2} - \bar\partial\varphi\wedge\partial\omega^{n-2})>0 \hspace{2ex} \mbox{and} \hspace{2ex} \sup\limits_X\varphi=0,\end{equation}

\noindent for a function $\varphi\,:\,X\rightarrow\R$, where $f$ is a given $C^{\infty}$ real-valued function. \\

\noindent $(a)$\, For any given $f$, are solutions $\varphi$ to $(\star)$ and (\ref{eqn:conditions}) unique?

\noindent $(b)$\, For any given $f$, let $\varphi$ be a $C^{\infty}$ solution of equation $(\star)$ subject to (\ref{eqn:conditions}). Are there uniform a priori $C^{\infty}$ estimates on $\varphi$ depending only on $(X,\,\omega,\, f)$?

\noindent $(c)$\, For any given $f$, does there exist a (unique) constant $c\in\R$ such that the equation

\begin{equation}\label{eqn:main-eq-const}\bigg[\bigg(\omega^{n-1} + i\partial\bar\partial\varphi\wedge\omega^{n-2} + \frac{i}{2}\,(\partial\varphi\wedge\bar\partial\omega^{n-2} - \bar\partial\varphi\wedge\partial\omega^{n-2})\bigg)^{\frac{1}{n-1}}\bigg]^n = e^{f+c}\,\omega^n\end{equation}

\noindent admits a $C^{\infty}$ solution $\varphi$ satisfying (\ref{eqn:conditions})? This solution is unique if the answer to $(a)$ is affirmative.

\end{Question}

 Note that in the special case of a K\"ahler metric $\omega$, $\partial\omega^{n-2}=0$ and $\bar\partial\omega^{n-2}=0$, so equation $(\star)$ simplifies to the equation

\begin{equation}\nonumber\bigg[\bigg( \omega^{n-1} + i\partial\bar\partial\varphi\wedge\omega^{n-2}\bigg)^{\frac{1}{n-1}}\bigg]^n = e^f\,\omega^n\hspace{1ex}(\star\star)\end{equation}

\noindent with initial conditions

\begin{equation}\label{eqn:conditions-second}\omega^{n-1} + i\partial\bar\partial\varphi\wedge\omega^{n-2} > 0 \hspace{2ex} \mbox{and} \hspace{2ex} \sup\limits_X\varphi=0.\end{equation}

\noindent Notice that for $n=2$, equation $(\star\star)$ is the classical Calabi-Yau equation. At the time of writing the first version of this paper, the author was unfortunately unaware of the works by Fu, Wang and Wu who had discussed in [FWW10a] and [FWW10b] the equation $(\star\star)$ and also unaware of the work of Tosatti and Weinkove who had completely solved equation $(\star\star)$ on compact K\"ahler manifolds in [TW13a]. However, the emphasis of the present work is firmly on the non-K\"ahler context and on the new equation $(\star)$ adapted to it. 

 Besides its applications to moduli spaces of Calabi-Yau $\partial\bar\partial$-manifolds outlined above, equation $(\star)$ would also contribute to the further study of the deformation properties of the Gauduchon and sG cones introduced and studied in  $\S$\,\ref{section:cones}. 

 In $\S$\,\ref{section:uniqueness} we prove the uniqueness of solutions to equation $(\star)$ subject to (\ref{eqn:conditions}). In $\S$\,\ref{section:linearisation} we calculate the linearisation of equation $(\star)$ and observe that its principal part is the Laplacian associated with a certain Hermitian metric on $X$. The following statement sums up these results (cf. Theorem \ref{The:uniqueness}, Proposition \ref{Prop:linearisation} and Corollary \ref{Cor:linearisation} for more precise wording).

\begin{The}\label{The:equation-summing} Fix a compact Hermitian manifold $(X,\,\omega)$, $\mbox{dim}_{\C}X=n\geq 2$. 

\vspace{1ex}

\noindent $(i)$\, Part $(a)$\, of Question \ref{Question:equation} has an affirmative answer. 

\vspace{1ex}

\noindent $(ii)$\, The principal part of the linearisation of equation $(\star)$ is 

\vspace{1ex}

\hspace{30ex} $\displaystyle\frac{(n-2)!}{n-1}\,\Delta_{\tilde\lambda},$

\vspace{1ex}

\noindent where $\Delta_{\tilde\lambda} = \mbox{tr}_{\tilde\lambda}(i\partial\bar\partial)$ is the Laplacian associated with the $C^{\infty}$ positive definite $(1,\,1)$-form $\tilde\lambda$ defined by the following relations\!\!:

\vspace{2ex}

\noindent $\rho:=\star\Lambda>0$, where $\Lambda:=\omega^{n-1} + i\partial\bar\partial\varphi\wedge\omega^{n-2} + \frac{i}{2}\,(\partial\varphi\wedge\bar\partial\omega^{n-2} - \bar\partial\varphi\wedge\partial\omega^{n-2})>0,$

\noindent $\lambda$ is the $(n-1)^{st}$ root of $(\Lambda_{\rho}\omega)\,\frac{\omega^{n-1}}{(n-1)!} - (\omega^n/\rho^n)^{\frac{1}{n-1}}\,\star(\star\rho)^{\frac{1}{n-1}}>0,$

\noindent $\tilde\lambda = \frac{1}{(n-1)!}\frac{\lambda}{\lambda^n/\omega^n}>0,$

\vspace{1ex}

\noindent where $\star = \star_{\omega}$ is the Hodge star operator associated with $\omega$.

\end{The}

 Since the principal part of the linearisation of equation $(\star)$ is a constant factor of a Laplacian, the local inversion theorem can be applied as in the case of the classical Calabi-Yau equation to prove the openness of the interval of solutions in the continuity method. The resemblance with the latter equation makes it likely for $(\star)$ to lend itself to a treatment through the standard techniques developed in the literature for the classical Monge-Amp\`ere equation in bidegree $(1,\, 1)$. We hope to be able to take up the study of parts $(b)$ and $(c)$ of Question \ref{Question:equation} in future work.

\vspace{2ex}

\noindent {\bf Acknowledgments.} The author is very grateful to Jean-Pierre Demailly from discussions with whom during the autumn of 2009 the idea of studying an equation of the Monge-Amp\`ere type in bidegree $(n-1,\,n-1)$ for geometric applications to the deformation theory first emerged. Many thanks are also due for the interest he has shown since then in two earlier forms of such an equation with which the author has experimented over several years before hitting upon the idea of considering equation $(\star)$ in the context of the Aeppli cohomology as best suited to the original objectives.   
 
 Some 10 days after this work had been posted on the arXiv, Valentino Tosatti and Ben Weinkove informed the author that they were about to post their preprint [TW13b] in which they made significant progress towards the resolution of equation $(\star)$. The author is very grateful to them for their work on this equation and for letting him know of the earlier works [FWW10a], [FWW10b] and [TW13a] of which he was unfortunately unaware at the time.

\section{Bott-Chern and Aeppli cohomologies}\label{section:A-BC} Let $(X,\,\omega)$ still denote a compact Hermitian manifold, $\mbox{dim}_{\C}X=n$. We will give a different interpretation of Proposition \ref{Prop:pluri+bal}. 

 The $4^{th}$ order Bott-Chern Laplacian $\Delta^{p, \, q}_{BC}:C^{\infty}_{p, \, q}(X, \, \C)\rightarrow C^{\infty}_{p, \, q}(X, \, \C)$ introduced by Kodaira and Spencer in [KS60, $\S.6$] (see also [Sch07, 2.c., p. 9-10]) as defined by

\begin{equation}\label{eqn:BC-Laplacian}\Delta^{p, \, q}_{BC}:=\partial^{\star}\partial + \bar\partial^{\star}\bar\partial + (\partial\bar\partial)^{\star}(\partial\bar\partial) + (\partial\bar\partial)(\partial\bar\partial)^{\star} + (\partial^{\star}\bar\partial)^{\star}(\partial^{\star}\bar\partial) + (\partial^{\star}\bar\partial)(\partial^{\star}\bar\partial)^{\star}\end{equation}

\noindent is elliptic and formally self-adjoint, so it induces a three-space decomposition

\begin{equation}\label{eqn:BC-3sp-decomp}C^{\infty}_{p, \, q}(X, \C)=\ker\Delta^{p, \, q}_{BC} \oplus \mbox{Im}\,\partial\bar\partial \oplus (\mbox{Im}\,\partial^{\star} + \mbox{Im}\,\bar\partial^{\star})\end{equation}

\noindent that is orthogonal w.r.t. the $L^2$ scalar product defined by $\omega$. We have 

\begin{equation}\label{eqn:BC-kernel}\ker\partial\cap\ker\bar\partial=\ker\Delta^{p, \, q}_{BC} \oplus \mbox{Im}\,\partial\bar\partial,\end{equation}

\noindent yielding the Hodge isomorphism $H^{p, \, q}_{BC}(X, \, \C)\simeq \ker\Delta^{p, \, q}_{BC}$. We also have

\begin{equation}\label{eqn:BC-image}\mbox{Im}\,\Delta^{p, \, q}_{BC} = \mbox{Im}\,\partial\bar\partial \oplus (\mbox{Im}\,\partial^{\star} + \mbox{Im}\,\bar\partial^{\star}).\end{equation}

 Similarly, the $4^{th}$ order Aeppli Laplacian $\Delta^{p, \, q}_A:C^{\infty}_{p, \, q}(X, \, \C)\rightarrow C^{\infty}_{p, \, q}(X, \, \C)$ (cf. [Sch07, 2.c., p. 9-10]) defined by

\begin{equation}\label{eqn:A-Laplacian}\Delta^{p, \, q}_A:=\partial\partial^{\star} + \bar\partial\bar\partial^{\star} + (\partial\bar\partial)^{\star}(\partial\bar\partial) + (\partial\bar\partial)(\partial\bar\partial)^{\star} + (\partial\bar\partial^{\star})(\partial\bar\partial^{\star})^{\star} + (\partial\bar\partial^{\star})^{\star}(\partial\bar\partial^{\star})\end{equation}

\noindent is elliptic and formally self-adjoint, so it induces a three-space decomposition

\begin{equation}\label{eqn:A-3sp-decomp}C^{\infty}_{p, \, q}(X, \C)=\ker\Delta^{p, \, q}_A \oplus (\mbox{Im}\partial + \mbox{Im}\bar\partial) \oplus \mbox{Im}(\partial\bar\partial)^{\star}\end{equation}

\noindent that is orthogonal w.r.t. the $L^2$ scalar product defined by $\omega$. We have  

\begin{equation}\label{eqn:A-kernel}\ker(\partial\bar\partial)=\ker\Delta^{p, \, q}_A \oplus (\mbox{Im}\,\partial + \mbox{Im}\,\bar\partial),\end{equation}

\noindent yielding the Hodge isomorphism $H^{p, \, q}_A(X, \, \C)\simeq \ker\Delta^{p, \, q}_A$. We also have

\begin{equation}\label{eqn:A-image}\mbox{Im}\,\Delta^{p, \, q}_A = (\mbox{Im}\partial + \mbox{Im}\bar\partial) \oplus \mbox{Im}(\partial\bar\partial)^{\star}.\end{equation}

 In what follows, ${\cal H}^{p,\,q}_{\Delta_{BC}}(X,\,\C):=\ker\Delta^{p,\,q}_{BC}\subset C^{\infty}_{p,\,q}(X,\,\C)$ will stand for the space of Bott-Chern-harmonic $(p,\,q)$-forms and ${\cal H}^{p,\,q}_{\Delta_A}(X,\,\C):=\ker\Delta^{p,\,q}_A\subset C^{\infty}_{p,\,q}(X,\,\C)$ will denote the space of Aeppli-harmonic $(p,\,q)$-forms, while the Laplacians will be simply written $\Delta_{BC}$ and $\Delta_A$ (without the superscripts) when no confusion is likely. The following statement sums up the basic properties of $H^{p,\,q}_{BC}(X,\,\C)$, $H^{p,\,q}_A(X,\,\C)$ and their harmonic counterparts, some of which already appear in [KS, $\S.6$] and in [Sch07].

\begin{The}\label{The:BC-A-3space} Let $(X,\,\omega)$ be a compact Hermitian manifold, $\mbox{dim}_{\C}X=n$.

\noindent $(i)$\, We have 

\vspace{1ex}

  $\ker(\partial\bar\partial)^{\star}={\cal H}^{p,\,q}_{\Delta_{BC}}(X,\,\C)\oplus(\mbox{Im}\,\partial^{\star} + \mbox{Im}\,\bar\partial^{\star})$ 

\noindent and 

$\ker\partial^{\star}\cap\ker\bar\partial^{\star} = {\cal H}^{p,\,q}_{\Delta_A}(X,\,\C)\oplus \mbox{Im}\,(\partial\bar\partial)^{\star}.$ \\

\noindent It follows that 

\vspace{1ex}

${\cal H}^{p,\,q}_{\Delta_{BC}}(X,\,\C) = \ker\partial\cap\ker\bar\partial\cap\ker(\partial\bar\partial)^{\star}$

\noindent and

${\cal H}^{p,\,q}_{\Delta_A}(X,\,\C) = \ker(\partial\bar\partial)\cap\ker\partial^{\star}\cap\ker\bar\partial^{\star}.$

\vspace{1ex}

\noindent In particular, for any $(p,\,q)$-form $\alpha$, the following equivalences hold\!\!:

\vspace{1ex}

$\Delta_{BC}\alpha=0 \iff \Delta_{BC}\bar{\alpha}=0 \hspace{2ex}\mbox{and}\hspace{2ex} \Delta_A\alpha=0 \iff \Delta_A\bar{\alpha}=0.$

\vspace{1ex}

\noindent (Note that $\overline{\Delta_{BC}}\neq\Delta_{BC}$ and $\overline{\Delta_A}\neq\Delta_A$ because of the last two terms in the definition of each of these Laplacians.)

\vspace{1ex}

\noindent $(ii)$\, Under the Hodge star isomorphism $\star=\star_{\omega}:C^{\infty}_{p,\,q}(X,\,\C)\rightarrow C^{\infty}_{n-q,\,n-p}(X,\,\C)$ defined by $\omega$, the Bott-Chern and Aeppli three-space decompositions (\ref{eqn:BC-3sp-decomp}) of $C^{\infty}_{p,\,q}(X,\,\C)$ and respectively (\ref{eqn:A-3sp-decomp}) of $C^{\infty}_{n-q,\,n-p}(X,\,\C)$ are related by the following three restrictions of $\star$ being isomorphisms\!\!:

\begin{equation}\label{eqn:star-BC-A-harm}\star : {\cal H}^{p,\,q}_{\Delta_{BC}} \stackrel{\simeq}{\longrightarrow} {\cal H}^{n-q,\,n-p}_{\Delta_A},\end{equation}

\begin{equation}\label{eqn:star-crossing-isom}\star : \mbox{Im}\,\partial\bar\partial \stackrel{\simeq}{\longrightarrow} \mbox{Im}\,(\partial\bar\partial)^{\star} \hspace{2ex}\mbox{and}\hspace{2ex} \star : (\mbox{Im}\,\partial^{\star} + \mbox{Im}\,\bar\partial^{\star})\stackrel{\simeq}{\longrightarrow} (\mbox{Im}\,\partial + \mbox{Im}\,\bar\partial).\end{equation}

\noindent Thus the resulting isomorphism in cohomology

\begin{equation}\label{eqn:BC-A-coh-isom}\star : H^{p,\,q}_{BC}(X,\,\C)\rightarrow H^{n-q,\,n-p}_A(X,\,\C)\end{equation} 

\noindent depends on the choice of the metric $\omega$.

\noindent $(iii)$ The following duality in cohomology

\begin{equation}\label{eqn:duality-BC-A-def}H^{p,\,q}_{BC}(X,\,\C)\times H^{n-p,\,n-q}_A(X,\,\C) \longrightarrow \C, \hspace{2ex}([\alpha]_{BC},\, [\beta]_A) \longmapsto \int\limits_X\alpha\wedge\beta\end{equation}

\noindent is well defined, canonical (i.e. independent of the metric $\omega$) and non-degenerate.

\end{The}

\noindent {\it Proof.} $(i)$\, The three-space decomposition (\ref{eqn:BC-3sp-decomp}) being orthogonal, we have

\vspace{1ex}

\hspace{7ex} ${\cal H}^{p,\,q}_{BC}(X,\,\C)\oplus(\mbox{Im}\partial^{\star} + \mbox{Im}\bar\partial^{\star}) = (\mbox{Im}\partial\bar\partial)^{\perp} = \ker(\partial\bar\partial)^{\star},$

\vspace{1ex}

\noindent where the last identity is standard. Similarly, the orthogonality of decomposition (\ref{eqn:A-3sp-decomp}) gives

\vspace{1ex}

\noindent ${\cal H}^{p,\,q}_A(X,\,\C)\oplus\mbox{Im}(\partial\bar\partial)^{\star} = (\mbox{Im}\partial + \mbox{Im}\bar\partial)^{\perp} = (\mbox{Im}\partial)^{\perp}\cap(\mbox{Im}\bar\partial)^{\perp} = \ker\partial^{\star}\cap\ker\bar\partial^{\star}.$

\vspace{1ex}

\noindent This proves the first two identities in part $(i)$. The remaining two identities in $(i)$ follow immediately from these and from (\ref{eqn:BC-kernel}) and respectively (\ref{eqn:A-kernel}).

\vspace{1ex}

\noindent $(ii)$\, The well-known identities $\partial^{\star} = -\star\bar\partial\star$ and $\bar\partial^{\star} = -\star\partial\star$ imply the inclusions

\vspace{1ex}

\hspace{4ex}$\star\bigg(\mbox{Im}\,(\partial\bar\partial)^{\star}\bigg)\subset\mbox{Im}\,\partial\bar\partial  \hspace{2ex} \mbox{and} \hspace{2ex} \star\bigg(\mbox{Im}\,\partial + \mbox{Im}\,\bar\partial\bigg)\subset\mbox{Im}\,\partial^{\star} + \mbox{Im}\,\bar\partial^{\star}$

\vspace{1ex}

\noindent and, combined with part $(i)$, they also imply the equivalences (cf. [Sch07])\!\!:

\begin{eqnarray}\nonumber u\in{\cal H}^{p,\,q}_{\Delta_{BC}} & \iff & \partial u=0,\, \bar\partial u=0,\, (\partial\bar\partial)^{\star}u=0 \\
\nonumber & \iff & \bar\partial^{\star}(\star u)=0,\, \partial^{\star}(\star u)=0,\, \partial\bar\partial (\star u)=0 \\
\nonumber & \iff & \star u\in{\cal H}^{n-q,\,n-p}_{\Delta_A}.\end{eqnarray}

\noindent This proves (\ref{eqn:star-BC-A-harm}), while (\ref{eqn:star-crossing-isom}) follows immediately using the above inclusions.

\noindent $(iii)$\, It is obvious that the metric $\omega$ does not feature in the definition of the pairing (\ref{eqn:duality-BC-A-def}). To show that the pairing (\ref{eqn:duality-BC-A-def}) is well defined, i.e. independent of the choice of representatives $\alpha$, $\beta$ of the respective Bott-Chern and Aeppli classes, let $\alpha\in C^{\infty}_{p,\,q}(X,\,\C)$ and $\beta\in C^{\infty}_{n-p,\,n-q}(X,\,\C)$ be such that $d\alpha=0$ and $\partial\bar\partial\beta=0$. Any representative of the Bott-Chern class $[\alpha]_{BC}$ is of the shape $\alpha + \partial\bar\partial\gamma$ for some $\gamma\in C^{\infty}_{p-1,\,q-1}(X,\,\C)$\!; we have

\vspace{1ex}

\hspace{10ex}$\displaystyle\int\limits_X(\alpha + \partial\bar\partial\gamma)\wedge\beta = \int\limits_X\alpha\wedge\beta - \int\limits_X\gamma\wedge\partial\bar\partial\beta = \int\limits_X\alpha\wedge\beta$

\vspace{1ex}

\noindent since $\partial\bar\partial\beta=0$. Similarly, any representative of the Aeppli class $[\beta]_A$ is of the shape $\beta + \partial u + \bar\partial v$ for some $u\in C^{\infty}_{n-p-1,\,n-q}(X,\,\C)$ and $v\in C^{\infty}_{n-p,\,n-q-1}(X,\,\C)$; we have

\vspace{1ex}

\hspace{1ex}$\displaystyle\int\limits_X\alpha\wedge(\beta + \partial u + \bar\partial v) = \int\limits_X\alpha\wedge\beta \pm \int\limits_X\partial\alpha\wedge u \pm \int\limits_X\bar\partial\alpha\wedge v = \int\limits_X\alpha\wedge\beta$

\vspace{1ex}

\noindent since $\partial\alpha=0$ and $\bar\partial\alpha=0$. 

 That the pairing (\ref{eqn:duality-BC-A-def}) is non-degenerate follows from the isomorphism (\ref{eqn:star-BC-A-harm}). Indeed, if $[\alpha]_{BC}\in H^{p,\,q}_{BC}(X,\,\C)$ is any class, let $\alpha$ denote its unique Bott-Chern-harmonic representative. Then $\Delta_A(\star\alpha)=0$ by (\ref{eqn:star-BC-A-harm}), hence $\Delta_A(\star\bar\alpha)=0$ by the last statement of part $(i)$, so $\star\bar\alpha$ is the unique Aeppli-harmonic representative of the class $[\star\bar\alpha]_A\in H_A^{n-p,\,n-q}(X,\,\C)$ and we have

\vspace{1ex}

\hspace{23ex}$\displaystyle\int\limits_X\alpha\wedge\star\bar\alpha = ||\alpha||^2 >0$

\vspace{1ex}

\noindent if $\alpha\neq 0$ (i.e. if $[\alpha]_{BC}\neq 0$ in $H^{p,\,q}_{BC}(X,\,\C)$). Similarly, if $[\beta]_A\in H_A^{n-p,\,n-q}(X,\,\C)$ is any class and $\beta$ denotes its Aeppli-harmonic representative, then $\Delta_{BC}(\star\bar\beta)=0$ by (\ref{eqn:star-BC-A-harm}) and the last statement of part $(i)$, while  $\int_X\beta\wedge\star\bar\beta=||\beta||^2>0$ if $\beta\neq 0$ (i.e. if $[\beta]_A\neq 0$). \hfill $\Box$

\vspace{2ex}

 We can now observe that for a pluriclosed metric, the balanced condition is equivalent to the Aeppli harmonicity.

\begin{Lem}\label{Lem:pluri-closed_A-harm} Let $\omega>0$ be a $C^{\infty}$ positive definite $(1,\,1)$-form on $X$ such that $\partial\bar\partial\omega=0$. The following equivalence holds\!\!:

$$\Delta_A\omega =0 \iff d\omega^{n-1}=0.$$

\end{Lem}

\noindent {\it Proof.} Since $\star\omega = \omega^{n-1}/(n-1)!$ and $d^{\star} = -\star d\star$, the balanced condition $d\omega^{n-1}=0$ is equivalent to $d^{\star}\omega=0$, hence to $\partial^{\star}\omega=0$ and $\bar\partial^{\star}\omega=0$. The contention is thus seen to follow from the vector space identity ${\cal H}^{1,\,1}_{\Delta_A}=\ker(\partial\bar\partial)\cap\ker\partial^{\star}\cap\ker\bar\partial^{\star}$ proved in part $(i)$ of Theorem \ref{The:BC-A-3space}.  \hfill $\Box$

 Thus Proposition \ref{Prop:pluri+bal} can be reworded in the following way.

\begin{Cor}\label{Cor:rewording-pluri+bal} Let $\omega>0$ be a Hermitian metric on $X$. Then

$$\omega \hspace{2ex} \mbox{is K\"ahler} \hspace{2ex} \iff \Delta_A\omega=0.$$

\end{Cor}

\section{Relations with the $\partial\bar\partial$-lemma}\label{section:ddbar-rel}

A $C^{\infty}$ positive definite $(1,\,1)$-form $\omega$ on $X$ is Hermitian-symplectic (cf. definition in [ST10]) iff there exists $\alpha\in C^{\infty}_{0,\,2}(X,\,\C)$ s.t. $d(\overline{\alpha^{0,\,2}} + \omega + \alpha^{0,\,2})=0$, which amounts to

\begin{equation}\label{eqn:HS-conditions}\exists\,\alpha^{0,\,2}\in C^{\infty}_{0,\,2}(X,\,\C) \hspace{1ex} \mbox{s.t.} \hspace{1ex} \partial\omega + \bar\partial\overline{\alpha^{0,\,2}}=0 \hspace{2ex} \mbox{and} \hspace{2ex} \partial\overline{\alpha^{0,\,2}}=0.\end{equation}

\noindent Indeed, in the real $3$-form $d(\overline{\alpha^{0,\,2}} + \omega + \alpha^{0,\,2})$ the components of types $(3,\,0)$ and $(0,\,3)$ are conjugate to each other and so are the components of types $(2,\,1)$ and $(1,\,2)$, so the vanishing of $d(\overline{\alpha^{0,\,2}} + \omega + \alpha^{0,\,2})$ is equivalent to the vanishing of its components of types $(2,\,1)$ and $(3,\,0)$.

 We now observe that on a $\partial\bar\partial$-manifold, the two conditions in (\ref{eqn:HS-conditions}) characterising the Hermitian-symplectic property reduce to the first one and that, consequently, the notions of Hermitian-symplectic and pluriclosed metrics coincide.

\begin{Lem}\label{Lem:HS=pluriclosed} Let $X$ be a compact $\partial\bar\partial$-manifold. For any Hermitian metric $\omega$, the following equivalences hold\!\!:

\begin{equation}\label{eqn:HS=pluriclosed}\nonumber\omega \hspace{1ex} \mbox{is Hermitian-symplectic} \stackrel{(a)}{\Longleftrightarrow} \partial\omega\in\mbox{Im}\,\bar\partial \stackrel{(b)}{\Longleftrightarrow} \partial\bar\partial\omega=0.\end{equation}

\end{Lem}

\noindent {\it Proof.} To prove the implication $\stackrel{(b)}{\Longleftarrow}$, suppose that $\partial\bar\partial\omega=0$, which means that $\partial\omega\in\ker\bar\partial$, hence $\partial\omega$ is a $d$-closed form of pure type $(2,\,1)$. Since $\partial\omega$ is $\partial$-exact, it must also be $\bar\partial$-exact by the $\partial\bar\partial$-assumption on $X$. The implication $\stackrel{(b)}{\Longrightarrow}$ is obvious.

 To prove the implication $\stackrel{(a)}{\Longleftarrow}$, suppose that $\partial\omega\in\mbox{Im}\,\bar\partial$ and let $\alpha^{2,\,0}\in C^{\infty}_{2,\,0}(X,\,\C)$ such that $\partial\omega=-\bar\partial\alpha^{2,\,0}$. Put $\alpha^{0,\,2}:=\overline{\alpha^{2,\,0}}$. Then $\partial\omega + \bar\partial\overline{\alpha^{0,\,2}}=0$. In view of (\ref{eqn:HS-conditions}), it remains to show that $\partial\alpha^{2,\,0}=0$. 

 Now $\partial\alpha^{2,\,0}$ is $\bar\partial$-closed since $\bar\partial(\partial\alpha^{2,\,0})=-\partial(\bar\partial\alpha^{2,\,0})=\partial^2\omega=0.$ Thus the $(3,\,0)$-form $\partial\alpha^{2,\,0}$ is $d$-closed and $\partial$-exact, hence it must also be $\bar\partial$-exact by the $\partial\bar\partial$-assumption on $X$. However, the only $\bar\partial$-exact $(3,\,0)$-form is zero, hence $\partial\alpha^{2,\,0}=0$.  

 The implication $\stackrel{(a)}{\Longrightarrow}$ is obvious in view of (\ref{eqn:HS-conditions}).  \hfill $\Box$

\vspace{3ex}

 It is well known that on any compact complex manifold $X$ and for any $(p,\,q)$, there are well-defined linear maps from the Bott-Chern cohomology group $H^{p,\,q}_{BC}(X,\,\C)$ to the Dolbeault, De Rham and Aeppli cohomology groups $H^{p,\,q}_{\bar\partial}(X,\,\C)$, $H^{p+q}_{DR}(X,\,\C)$ and resp. $H^{p,\,q}_A(X,\,\C)$\!\!: 

$$[\alpha]_{BC}\longmapsto [\alpha]_{\bar\partial}, \hspace{2ex} [\alpha]_{BC}\longmapsto \{\alpha\}, \hspace{2ex} [\alpha]_{BC}\longmapsto [\alpha]_A,$$

\noindent and a well-defined linear map from the Dolbeault to the Aeppli cohomology\!:

\vspace{1ex}

\hspace{6ex} $H^{p,\,q}_{\bar\partial}(X,\,\C)\longrightarrow H^{p,\,q}_A(X,\,\C), \hspace{2ex} [\alpha]_{\bar\partial}\longmapsto [\alpha]_A.$

\vspace{1ex}

\noindent These maps are neither injective, nor surjective in general. However, if the $\partial\bar\partial$-lemma holds on $X$, the map to De Rham cohomology is injective while the others are isomorphisms. In the same vein, still denoting De Rham classes by $\{\,\,\,\}$, we have the following.

\begin{The}\label{The:ddbar-Aeppli} Let $X$ be a compact $\partial\bar\partial$-manifold, $\mbox{dim}_{\C}X=n$. 

\vspace{1ex}

\noindent $(a)$\, Every Aeppli cohomology class contains a $d$-closed representative. 

\noindent $(b)$\, For any $p, q=0, 1, \dots , n$, there is a canonical injective linear map\!\!:

\begin{equation}\label{eqn:A-DR-inj} T^{p,\,q}:H^{p,\,q}_A(X,\,\C)\hookrightarrow H^{p+q}_{DR}(X,\,\C), \hspace{2ex} [\alpha]_A\longmapsto \{\alpha\},\end{equation}

\noindent where $\alpha$ is any $d$-closed $(p,\,q)$-form representing the Aeppli class $[\alpha]_A$ whose existence is guaranteed by $(a)$.

\noindent $(c)$\, For any $k=0, 1, \dots , 2n$, there is a canonical isomorphism\!\!:

\begin{eqnarray}\label{eqn:Aeppli-Hodge-decomp}H^k_{DR}(X,\,\C) & \simeq & \bigoplus\limits_{p+q=k}H^{p,\,q}_A(X,\,\C)\\
\nonumber      \bigg\{\sum\limits_{p+q=k}\alpha^{p,\,q}\bigg\} & \rotatebox{180}{$\longmapsto$} & \sum\limits_{p+q=k}[\alpha^{p,\,q}]_A,\end{eqnarray}

\noindent where each $\alpha^{p,\,q}$ is a $d$-closed representative of the Aeppli class $[\alpha^{p,\,q}]_A$, that can well be called the {\bf Hodge-Aeppli decomposition}. Note that the Aeppli cohomology analogue of the Hodge symmetry always (even without the $\partial\bar\partial$-assumption on $X$) holds trivially, i.e. $H_A^{p,\,q}(X,\,\C)=\overline{H_A^{q,\,p}(X,\,\C)}$ for all $p, q$.

\end{The}

\noindent {\it Proof.} $(a)$\, Let $\alpha$ be a $(p,\,q)$-form such that $\partial\bar\partial\alpha=0$. We have to prove the existence of a $(p-1,\,q)$-form $\beta$ and of a $(p,\,q-1)$-form $\gamma$ such that $d(\alpha + \partial\beta + \bar\partial\gamma)=0$. The last identity translates to

$$\partial\alpha = -\partial\bar\partial\gamma  \hspace{2ex} \mbox{and} \hspace{2ex} \bar\partial\alpha = \partial\bar\partial\beta.$$

\noindent We are thus reduced to showing that $\partial\alpha$ and $\bar\partial\alpha$ are $\partial\bar\partial$-exact. Both $\partial\alpha$ and $\bar\partial\alpha$ are of pure types ($(p+1,\,q)$, resp. $(p,\,q+1)$) and $d$-closed (thanks to the assumption $\partial\bar\partial\alpha=0$), while $\partial\alpha$ is $\partial$-exact and $\bar\partial\alpha$ is $\bar\partial$-exact, so both must be $\partial\bar\partial$-exact by the $\partial\bar\partial$-lemma that holds on $X$ by hypothesis.  

\vspace{1ex}

\noindent $(b)$\, First, we have to show that $T^{p,\, q}$ is independent of the choice of $d$-closed representative of the Aeppli class $[\alpha]_A$. Let $\tilde{\alpha}, \tilde{\beta}\in C^{\infty}_{p,\,q}(X,\,\C)$ be $d$-closed forms representing the same Aeppli class, i.e.

$$d\tilde{\alpha}=d\tilde{\beta}=0 \hspace{2ex} \mbox{and} \hspace{2ex} \tilde{\alpha}-\tilde{\beta}=\partial u + \bar\partial v,$$

\noindent for some $(p-1,\,q)$-form $u$ and some $(p,\,q-1)$-form $v$. It follows that

\vspace{1ex}

$0=\partial(\tilde{\alpha}-\tilde{\beta})=\partial(\bar\partial v)$, hence $\bar\partial v$ is $\partial$-closed, hence $\bar\partial v$ is $d$-closed.

\vspace{1ex}

\noindent Since $\bar\partial v$ is obviously a $\bar\partial$-exact pure-type form, the $\partial\bar\partial$-assumption on $X$ implies that $\bar\partial v\in \mbox{Im}\,\partial\bar\partial$. Similarly, we have

\vspace{1ex}

$0=\bar\partial(\tilde{\alpha}-\tilde{\beta})=\bar\partial(\partial u)$, hence $\partial u$ is $\bar\partial$-closed, hence $\partial u$ is $d$-closed.

\vspace{1ex}

\noindent Since $\partial u$ is obviously a $\partial$-exact pure-type form, the $\partial\bar\partial$-assumption on $X$ implies that $\partial u\in \mbox{Im}\,\partial\bar\partial$.

 Putting together the last two pieces of information, we find that

\vspace{1ex}

\hspace{6ex} $\tilde{\alpha} - \tilde{\beta} = \partial u + \bar\partial v \in\mbox{Im}\,\partial\bar\partial\subset\mbox{Im}\,d.$

\vspace{1ex}

\noindent Thus $\tilde{\alpha}$ and $\tilde{\beta}$ are $d$-cohomologous, so they define the same De Rham cohomology class $\{\tilde\alpha\}=\{\tilde\beta\}\in H^{p+q}_{DR}(X,\,\C)$, i.e. $T^{p,\,q}([\tilde\alpha]_A) = T^{p,\,q}([\tilde\beta]_A)$.

It remains to show that $T^{p,\,q}$ is injective. Let $\alpha\in C^{\infty}_{p,\,q}(X,\,\C)$ such that $d\alpha=0$ and $T^{p,\,q}([\alpha]_A)=\{\alpha\}=0$. The last identity means that $\alpha$ is $d$-exact. By the $\partial\bar\partial$-assumption on $X$, $\alpha$ must also be $\partial\bar\partial$-exact. In particular, $\alpha\in\mbox{Im}\,\partial + \mbox{Im}\,\bar\partial$, which means that $[\alpha]_A=0$.

\vspace{1ex}

\noindent $(c)$\, If $T : \bigoplus\limits_{p+q=k}H^{p,\,q}_A(X,\,\C)\longrightarrow H^k(X,\,\C)$ is the linear map $T=\sum\limits_{p+q=k}T^{p,\,q}$, then $T$ is injective since each $T^{p,\,q}$ is and the images in $H^k(X,\,\C)$ of any two different $H^{p,\,q}_A(X,\,\C)$ meet only at zero. Since $X$ is compact and $\Delta:=dd^{\star} + d^{\star}d$ and $\Delta_A$ (defined for any Hermitian metric on $X$) are elliptic, all the vector spaces involved are finite-dimensional, so the injectivity of $T$ implies

\begin{equation}\label{eqn:dim-eq-A-DR}\sum\limits_{p+q=k}\mbox{dim}_{\C}H^{p,\,q}_A(X,\,\C)\leq\mbox{dim}_{\C}H^k(X,\,\C).\end{equation}

\noindent On the other hand, the $\partial\bar\partial$-assumption on $X$ implies that $H^{p,\,q}_A(X,\,\C)$ is isomorphic to $H^{p,\,q}_{\bar\partial}(X,\,\C)$ for every $p, q$ and that 

$$\sum\limits_{p+q=k}\mbox{dim}_{\C}H^{p,\,q}_{\bar\partial}(X,\,\C)=\mbox{dim}_{\C}H^k(X,\,\C).$$

\noindent Thus equality holds in (\ref{eqn:dim-eq-A-DR}) for all $k$, hence the injective map $T=\sum\limits_{p+q=k}T^{p,\,q}$ must be an isomorphism. \hfill $\Box$

\vspace{3ex}

 For any compact complex manifold $X$ (not necessarily $\partial\bar\partial$) and any $p$, the space $H_A^{p,\,p}(X,\,\C)$ is stable under conjugation, so we can define $H_A^{p,\,p}(X,\,\R)\subset H_A^{p,\,p}(X,\,\C)$ to be the real subspace of real Aeppli $(p,\,p)$-classes (i.e. classes such that $\overline{[u]}_A=[u]_A$). Thus $H_A^{p,\,p}(X,\,\R)$ is the subspace of classes representable by a real $(p,\,p)$-form. Note that thanks to the last statement in $(i)$ of Theorem \ref{The:BC-A-3space}, for any Hermitian metric on $X$, the Aeppli-harmonic representative of a real Aeppli $(p,\,p)$-class is real.

\section{Resolution of the $\partial\bar\partial$ equation}\label{section:ddbar-equation}

 Let $(X,\,\omega)$ be a compact Hermitian manifold. The Bott-Chern Laplacian $\Delta_{BC}$ can be used to derive an explicit formula and an estimate for the minimal $L^2$-norm solution of the $\partial\bar\partial$-equation on $X$ that parallels standard formulae for the minimal solutions of the $d$, $\partial$ and $\bar\partial$-equations known in terms of $\Delta$, $\Delta'$ and resp. $\Delta''$. Similar uses of $\Delta_{BC}$ have been made in [KS60] and [FLY,$\S.4$].

 It will prove useful later on to consider as well the following $4^{th}$ order {\it real} Laplace-type operator that we will call {\it $\partial\bar\partial$-Laplacian}\!\!:

\begin{equation}\label{eqn:KS-Laplacian}\Delta_{\partial\bar\partial}:= (\partial\bar\partial)(\partial\bar\partial)^{\star} + (\partial\bar\partial)^{\star}(\partial\bar\partial).\end{equation}

\noindent It is obvious that $\overline{\Delta_{\partial\bar\partial}} = \Delta_{\partial\bar\partial}$ and that

\begin{equation}\label{eqn:ker-ddbar-Laplacian}\ker\Delta_{\partial\bar\partial}=\ker(\partial\bar\partial)\cap\ker(\partial\bar\partial)^{\star}\supset\ker\Delta_{BC}=\ker\partial\cap\ker\bar\partial\cap\ker(\partial\bar\partial)^{\star}.\end{equation}

\begin{The}\label{The:min-sol-ddbar} Fix a compact Hermitian manifold $(X,\,\omega)$. For any $C^{\infty}$ $(p,\,q)$-form $v\in\mbox{Im}\,(\partial\bar\partial)$, the (unique) minimal $L^2$-norm solution of the equation

\begin{equation}\label{eqn:ddbar-eq}\partial\bar\partial u=v\end{equation}

\noindent is given by the formula

\begin{equation}\label{eqn:ddbar-eq-solformula1}u = (\partial\bar\partial)^{\star}\Delta_{BC}^{-1}v,\end{equation}

\noindent as well as by the formula

\begin{equation}\label{eqn:ddbar-eq-solformula2}u = (\partial\bar\partial)^{\star}\Delta_{\partial\bar\partial}^{-1}v,\end{equation}

\noindent while its $L^2$-norm is estimated as

\begin{equation}\label{eqn:ddbar-eq-sol-est}||u||^2\leq\frac{1}{\lambda}\,||v||^2,\end{equation}

\noindent where $\Delta_{BC}^{-1}$ (resp. $\Delta_{\partial\bar\partial}^{-1}$) denotes the Green operator of $\Delta_{BC}$ (resp. of $\Delta_{\partial\bar\partial}$) and $\lambda>0$ is the smallest positive eigenvalue of $\Delta_{BC}$. Furthermore, we have

\begin{equation}\label{eqn:ddbar-eq-additional}\partial\Delta_{BC}^{-1}v=0 \hspace{2ex} \mbox{and} \hspace{2ex} \bar\partial\Delta_{BC}^{-1}v=0.\end{equation}

\end{The}

\noindent {\it Proof.} Let $w:=\Delta_{BC}^{-1}v$, i.e. $w$ is the unique $(p,\,q)$-form characterised by the following two properties

\begin{equation}\label{eqn:w-characterisation}\Delta_{BC}w=v  \hspace{2ex} \mbox{and} \hspace{2ex} w\perp\ker\Delta_{BC}.\end{equation}

\noindent By the definition (\ref{eqn:BC-Laplacian}) of $\Delta_{BC}$, the identity $\Delta_{BC}w=v=\partial\bar\partial u$ is equivalent to

\vspace{1ex}

\hspace{20ex} $A_1 + (A_2 + A_3) = 0,$ where

\begin{eqnarray}\nonumber A_1 & := & \partial\bar\partial\bigg((\partial\bar\partial)^{\star}w - u\bigg)\in\mbox{Im}\,\partial\cap\mbox{Im}\,\bar\partial,\\ 
\nonumber A_2 & := & \partial^{\star}\partial w + (\partial^{\star}\bar\partial)(\partial^{\star}\bar\partial)^{\star}w\in\mbox{Im}\,\partial^{\star},\\
\nonumber A_3 & := & \bar\partial^{\star}\bar\partial w + (\partial\bar\partial)^{\star}(\partial\bar\partial)w + (\partial^{\star}\bar\partial)^{\star}(\partial^{\star}\bar\partial)\in\mbox{Im}\,\bar\partial^{\star}.\end{eqnarray}

\noindent Since $\mbox{Im}\,\partial\perp\mbox{Im}\,\partial^{\star}$ and $\mbox{Im}\,\bar\partial\perp\mbox{Im}\,\bar\partial^{\star}$, we infer that $A_1\perp A_2$ and $A_1\perp A_3$, hence $A_1\perp(A_2+A_3)$. It follows that the identity $\Delta_{BC}w=v=\partial\bar\partial u$ is equivalent to $A_1=0$ and $A_2 + A_3=0$. Note that $A_1=0$ amounts to

\begin{equation}\label{eqn:A1vanishing}(\partial\bar\partial)^{\star}w - u \in\ker(\partial\bar\partial).\end{equation}

\noindent Meanwhile, the solutions of equation (\ref{eqn:ddbar-eq}) are unique up to $\ker(\partial\bar\partial)$, so if $u$ is the minimal $L^2$-norm solution, then $u\in\ker(\partial\bar\partial)^{\perp}=\mbox{Im}\,(\partial\bar\partial)^{\star}$. Thus

\begin{equation}\label{eqn:perp}(\partial\bar\partial)^{\star}w - u\in\mbox{Im}\,(\partial\bar\partial)^{\star}.\end{equation}

\noindent Now, $\ker(\partial\bar\partial)$ and $\mbox{Im}\,(\partial\bar\partial)^{\star}$ are mutually orthogonal, so thanks to (\ref{eqn:A1vanishing}) and (\ref{eqn:perp}), the identity $A_1=0$ is equivalent to $(\partial\bar\partial)^{\star}w - u=0$. This proves formula (\ref{eqn:ddbar-eq-solformula1}). On the other hand, the identity $A_2+A_3=0$ implies $\langle\langle A_2 + A_3,\,w\rangle\rangle=0$ which translates to

\begin{equation}\nonumber ||\partial w||^2 + ||\bar\partial^{\star}\partial w||^2 + ||\bar\partial w||^2 + ||\partial\bar\partial w||^2 + ||\partial^{\star}\bar\partial w||^2=0.\end{equation}

\noindent This amounts to $\partial w=0$ and $\bar\partial w=0$, proving (\ref{eqn:ddbar-eq-additional}).

 Let us now estimate the $L^2$ norm of $u=(\partial\bar\partial)^{\star}\Delta_{BC}^{-1}v$. We have

\begin{eqnarray}\nonumber||u||^2 & = & \langle\langle(\partial\bar\partial)(\partial\bar\partial)^{\star}\Delta_{BC}^{-1}v,\, \Delta_{BC}^{-1}v\rangle\rangle \stackrel{(a)}{=} \langle\langle\Delta_{BC}\Delta_{BC}^{-1}v,\,\Delta_{BC}^{-1}v\rangle\rangle \\
\nonumber & = & \langle\langle v,\,\Delta_{BC}^{-1}v\rangle\rangle\stackrel{(b)}{\leq}\frac{1}{\lambda}\,||v||^2,\end{eqnarray}

\noindent where identity $(a)$ follows from (\ref{eqn:BC-Laplacian}) and from the identities

\begin{eqnarray}\nonumber \partial^{\star}\partial\Delta_{BC}^{-1}v & = & 0,  \hspace{2ex} (\partial^{\star}\bar\partial)(\partial^{\star}\bar\partial)^{\star}\Delta_{BC}^{-1}v=0,\\
\nonumber \bar\partial^{\star}\bar\partial\Delta_{BC}^{-1}v & = & 0, \hspace{2ex} (\partial\bar\partial)^{\star}(\partial\bar\partial)\Delta_{BC}^{-1}v = 0,  \hspace{2ex} (\partial^{\star}\bar\partial)^{\star}(\partial^{\star}\bar\partial)\Delta_{BC}^{-1}v = 0,\end{eqnarray}

\noindent all of which are consequences of $\partial\Delta_{BC}^{-1}v = 0$ and of $\bar\partial\Delta_{BC}^{-1}v = 0$ already proved as (\ref{eqn:ddbar-eq-additional}). Inequality $(b)$ follows from $v\perp\ker\Delta_{BC}$ since $v\in\mbox{Im}\,(\partial\bar\partial)\subset\mbox{Im}\,\Delta_{BC}$ (see (\ref{eqn:BC-image})). Estimate (\ref{eqn:ddbar-eq-sol-est}) is proved.

 It remains to prove formula (\ref{eqn:ddbar-eq-solformula2}). The minimal $L^2$-norm solution $u$ of equation (\ref{eqn:ddbar-eq}) is the unique $(p-1,\, q-1)$-form $u$ satisfying the following two properties

\begin{equation}\label{eqn:prop-min-sol}\partial\bar\partial u=v \hspace{2ex}\mbox{and}\hspace{2ex} u\in\ker(\partial\bar\partial)^{\perp} = \mbox{Im}\,(\partial\bar\partial)^{\star}.\end{equation}

\noindent Let $u':=(\partial\bar\partial)^{\star}\Delta_{\partial\bar\partial}^{-1}v$. To prove that $u=u'$, we have to prove that $u'$ satisfies the two properties of (\ref{eqn:prop-min-sol}). Since it obviously satisfies the latter property, we are reduced to showing that $\partial\bar\partial u'=v$. We have

\noindent $\partial\bar\partial u'= (\partial\bar\partial)(\partial\bar\partial)^{\star}\Delta_{\partial\bar\partial}^{-1}v \stackrel{(i)}{=} \bigg((\partial\bar\partial)(\partial\bar\partial)^{\star} + (\partial\bar\partial)^{\star}(\partial\bar\partial)\bigg)\Delta_{\partial\bar\partial}^{-1}v = \Delta_{\partial\bar\partial}\Delta_{\partial\bar\partial}^{-1}v = v,$

\noindent where identity $(i)$ above follows from the commutation of $\partial\bar\partial$ with $\Delta_{\partial\bar\partial}$\!\!:

$$(\partial\bar\partial)\Delta_{\partial\bar\partial} = (\partial\bar\partial)(\partial\bar\partial)^{\star}(\partial\bar\partial) = \Delta_{\partial\bar\partial}(\partial\bar\partial),$$

\noindent which implies that $\partial\bar\partial$ and $\Delta_{\partial\bar\partial}^{-1}$ commute, which in turn implies the following identities

$$(\partial\bar\partial)^{\star}(\partial\bar\partial)\Delta_{\partial\bar\partial}^{-1}v = (\partial\bar\partial)^{\star}\Delta_{\partial\bar\partial}^{-1}(\partial\bar\partial v) =0$$

\noindent since $\partial\bar\partial v=0$ by assumption ($v$ is even assumed $\partial\bar\partial$-exact.)   \hfill $\Box$

\section{Cones of classes of metrics}\label{section:cones}

 Let $X$ be a compact complex manifold ($\mbox{dim}_{\C}X=n$). The canonical map

\begin{equation}\label{sG-map-def}T:H_A^{n-1,\,n-1}(X,\,\C)\longrightarrow H_{\bar\partial}^{n,\,n-1}(X,\,\C), \hspace{2ex} [\Omega]_A \mapsto [\partial\Omega]_{\bar\partial}\end{equation}

\noindent is well defined. Indeed, if $\Omega\in C^{\infty}_{n-1,\,n-1}(X,\,\C)$ defines an Aeppli cohomology class, then $\partial\bar\partial\Omega=0$, which amounts to $\partial\Omega$ being $\bar\partial$-closed, hence $\partial\Omega$ defines a Dolbeault cohomology class of bidegree $(n,\,n-1)$. If $\Omega_1$, $\Omega_2$ are two representatives of the same $(n-1,\,n-1)$ Aeppli class, then $\Omega_1 = \Omega_2 + \partial u + \bar\partial v$ for some forms $u, v$ of types $(n-2,\,n-1)$, resp. $(n-1,\,n-2)$. Thus $\partial\Omega_1 = \partial\Omega_2 + \bar\partial(- \partial v)$, hence $\partial\Omega_1$ and $\partial\Omega_2$ represent the same Dolbeault cohomology class, showing that $T([\Omega]_A)$ does not depend on the choice of representative of the Aeppli class $[\Omega]_A$.

 Now let $\omega$ be a Gauduchon metric on $X$. Then $\partial\bar\partial\omega^{n-1}=0$, so $\omega^{n-1}$ defines an Aeppli cohomology class $[\omega^{n-1}]_A\in H_A^{n-1,\,n-1}(X,\,\R)$ that will be called the {\bf Aeppli-Gauduchon class} associated with $\omega$. It is clear that

$$[\omega^{n-1}]_A\in\ker T \iff \partial\omega^{n-1}\in\mbox{Im}\,\bar\partial \iff \omega \hspace{1ex}\mbox{is a strongly Gauduchon metric},$$

\noindent the last equivalence being precisely the definition of a strongly Gauduchon (sG) metric (cf. [Pop09]). This shows that the strongly Gauduchon property is cohomological in the sense that either all Gauduchon metrics $\omega$ with $\omega^{n-1}$ lying in a given Aeppli class are strongly Gauduchon, or none of them is.

\begin{Def}\label{Def:sG-classes} $(i)$\, An {\bf sG class} on $X$ is an Aeppli-Gauduchon class lying in $\ker T$, i.e. any Aeppli cohomology class $[\omega^{n-1}]_A\in H^{n-1,\,n-1}_A(X,\,\R)$ representable by the $(n-1)^{st}$ power of a strongly Gauduchon metric $\omega$. 

\noindent $(ii)$\, The {\bf Gauduchon cone} of $X$ is the set ${\cal G}_X\subset H^{n-1,\,n-1}_A(X,\,\R)$ of Aeppli-Gauduchon classes, i.e. the convex cone of Aeppli classes $[\omega^{n-1}]_A$ of $(n-1)^{st}$ powers $\omega^{n-1}$ of Gauduchon metrics $\omega$.

\noindent $(iii)$\, The {\bf sG cone} of $X$ is the set ${\cal SG}_X\subset H^{n-1,\,n-1}_A(X,\,\R)$ of sG classes, i.e. the subcone of the Gauduchon cone defined as the intersection 

$${\cal SG}_X =  {\cal G}_X\cap\ker T\subset{\cal G}_X\subset H^{n-1,\,n-1}_A(X,\,\R).$$

\end{Def}

 Note that the subsets of $H^{n-1,\,n-1}_A(X,\,\R)$ defined above are indeed convex cones as follows by taking $(n-1)^{st}$ roots. For example, if $[\omega_1^{n-1}]_A, [\omega_2^{n-1}]_A\in{\cal G}_X$, then $[\omega_1^{n-1}]_A + [\omega_2^{n-1}]_A=[\omega^{n-1}]_A\in{\cal G}_X$ where $\omega>0$ is the $(n-1)^{st}$ root of $\omega_1^{n-1}+\omega_2^{n-1}>0$.

 We easily infer the following.

\begin{Obs}\label{Obs:Gcone-open}The Gauduchon cone ${\cal G}_X$ is an open subset of $H^{n-1,\,n-1}_A(X,\,\R)$. 

\end{Obs}

\noindent {\it Proof.} Let us equip the finite-dimensional vector space $H^{n-1,\,n-1}_A(X,\,\R)$ with an arbitrary norm $||\hspace{1ex}||$ (e.g. the Euclidian norm after we have fixed a basis; at any rate, all the norms are equivalent). Let $[\omega^{n-1}]_A\in{\cal G}_X$ be an arbitrary element, where $\omega>0$ is some Gauduchon metric on $X$. Let $\alpha\in H^{n-1,\,n-1}_A(X,\,\R)$ be a class such that $||\alpha - [\omega^{n-1}]_A||<\varepsilon$ for some small $\varepsilon>0$. Fix any Hermitian metric $\omega_0$ on $X$ and consider the Aeppli Laplacian $\Delta_A$ defined by $\omega_0$ inducing the Hodge isomorphism $H^{n-1,\,n-1}_A(X,\,\R)\simeq{\cal H}^{n-1,\,n-1}_{\Delta_A}(X,\,\R)$. Let $\Omega_{\alpha}\in{\cal H}^{n-1,\,n-1}_{\Delta_A}(X,\,\R)$ be the $\Delta_A$-harmonic representative of the class $\alpha$. Since $\omega^{n-1}\in\ker(\partial\bar\partial)$, (\ref{eqn:A-kernel}) gives a unique decomposition

$$\omega^{n-1} = \Omega + (\partial u + \bar\partial v) \hspace{2ex} \mbox{with}\hspace{1ex} \Delta_A\Omega=0.$$

\noindent If we set $\Gamma:=\Omega_{\alpha} + (\partial u + \bar\partial v)$ (with the same forms $u,v$ as for $\omega^{n-1}$), then $\partial\bar\partial\Gamma=0$, $\Gamma$ represents the Aeppli class $\alpha$ and we have

$$||\Gamma - \omega^{n-1}||_{C^0} = ||\Omega_{\alpha} - \Omega||_{C^0}\leq C\,||\alpha - [\omega^{n-1}]_A|| < C\varepsilon,$$

\noindent for some constant $C>0$ induced by the Hodge isomorphism. (We have chosen the $C^0$ norm on ${\cal H}^{n-1,\,n-1}_{\Delta_A}(X,\,\R)$ only for the sake of convenience.) Thus, if $\varepsilon>0$ is chosen sufficiently small, the $(n-1,\,n-1)$-form $\Gamma$ must be positive definite since $\omega^{n-1}$ is, so there exists a unique positive definite $(1,\,1)$-form $\gamma$ such that $\gamma^{n-1}=\Gamma$. Thus $\gamma$ is a Gauduchon metric and $\gamma^{n-1}$ represents the original Aeppli class $\alpha$, so $\alpha\in{\cal G}_X$. \hfill $\Box$

\vspace{3ex}

 Note that the Gauduchon cone is never empty since Gauduchon metrics exist on any compact complex manifold $X$ (cf. [Gau77]), while the sG cone of $X$ is empty if and only if $X$ is not an sG manifold. On the other hand, the sG cone of any $\partial\bar\partial$-manifold $X$ is maximal, i.e. ${\cal SG}_X = {\cal G}_X$, since on a $\partial\bar\partial$-manifold every Gauduchon metric is strongly Gauduchon (cf. [Pop09]). So we have the following implications\!\!:

\vspace{1ex}

\hspace{3ex} $X\,\,\mbox{is a}\,\,\partial\bar\partial\mbox{-manifold}\implies {\cal SG}_X = {\cal G}_X \implies X\,\,\mbox{is an sG-manifold}.$

\vspace{1ex}

\noindent In our opinion, compact complex manifolds $X$ for which ${\cal SG}_X = {\cal G}_X$ deserve further study. For example, their behaviour under deformations of the complex structure warrants being understood.

\begin{Obs}\label{Obs:cones-equality} The equality of cones ${\cal SG}_X = {\cal G}_X$ is equivalent to the following very special case of the $\partial\bar\partial$-lemma\!: every smooth $d$-closed $\partial$-exact $(n,\,n-1)$-form on $X$ is $\bar\partial$-exact (i.e. $T\equiv 0$).

\end{Obs}

\noindent {\it Proof.} Since ${\cal SG}_X = {\cal G}_X\cap\ker T = {\cal G}_X\cap(\ker T\cap H_A^{n-1,\,n-1}(X,\,\R))$ and ${\cal G}_X$ is open in $H_A^{n-1,\,n-1}(X,\,\R)$, the equality ${\cal SG}_X = {\cal G}_X$ is equivalent to $\ker T\cap H_A^{n-1,\,n-1}(X,\,\R) = H_A^{n-1,\,n-1}(X,\,\R)$, i.e. to $H_A^{n-1,\,n-1}(X,\,\R)\subset\ker T$. Since $\ker T$ is a $\C$ vector subspace of $H_A^{n-1,\,n-1}(X,\,\C)$, the last inclusion amounts to $\ker T = H_A^{n-1,\,n-1}(X,\,\C)$, i.e. to $T$ being identically zero.  \hfill $\Box$

\vspace{2ex}

It is worth noticing that there are examples of compact complex manifolds $X$ whose Gauduchon cone is the whole space $H^{n-1,\,n-1}_A(X,\,\R)$. In this case, we will say that the Gauduchon cone {\bf degenerates}. If $X$ is the connected sum $\sharp_k(S^3\times S^3)$ of $k\geq 2$ copies of $S^3\times S^3$, it was shown in [FLY12, Corollary 1.3] that the complex structure constructed on $X$ in [Fri91] and [LT96] by ``conifold transitions'' admits a balanced metric $\omega$. Since $\mbox{dim}_{\C}X=3$, $\omega^2$ defines a De Rham cohomology class in $H^4(X,\,\C)$. However, $H^4(X,\,\C)=0$ for this particular $X$, so $\omega^2$ must be $d$-exact. In particular, $\omega^2\in\mbox{Im}\,\partial + \mbox{Im}\,\bar\partial$, hence $[\omega^2]_A=0$. Since $\omega$ is necessarily a Gauduchon metric on $X$, it follows that ${\cal G}_X$ contains the origin, hence due to being open it must contain a neighbourhood of $0$ in $H^{2,\,2}_A(X,\,\R)$. Then ${\cal G}_X = H^{2,\,2}_A(X,\,\R)$ by the convex cone property of ${\cal G}_X$. It would be interesting to know whether the identity ${\cal G}_X=H^{n-1,\,n-1}_A(X,\,\R)$ (which is clearly equivalent to $0\in{\cal G}_X$ by the above arguments) can hold when $H^2(X,\,\C)\neq 0$ or $H^{n-1,\,n-1}_A(X,\,\R)\neq 0$.

 The following statement shows that the manifolds whose Gauduchon cone degenerates are rather exotic.

\begin{Prop}\label{Prop:exact-bal} Let $X$ be a compact complex manifold, $\mbox{dim}_{\C}X=n$. 

\vspace{1ex}

\noindent $(a)$\,The following three statements are equivalent.

\vspace{1ex}

 $(i)$\, There exists a $d$-exact $C^{\infty}$ $(n-1,\,n-1)$-form $\Omega>0$ on $X$ (henceforth called a {\bf degenerate balanced structure}).

\vspace{1ex}

 $(ii)$\, There exists no nonzero $d$-closed $(1,\,1)$-current $T\geq 0$ on $X$.

\vspace{1ex}

 $(iii)$\, The Gauduchon cone of $X$ degenerates\!\!: ${\cal G}_X= H_A^{n-1,\, n-1}(X,\,\R)$.

\vspace{1ex}

\noindent Furthermore, if any of the above three equivalent properties holds, $X$ cannot be a class ${\cal C}$ manifold.

\vspace{1ex}

\noindent $(b)$\, If $H^2(X,\,\C)=0$, the following equivalence holds

\vspace{1ex}

\hspace{6ex}$X\,\, \mbox{is an sG manifold} \hspace{1ex}\iff\hspace{1ex} X\,\, \mbox{is a balanced manifold}$

\vspace{1ex}

 \noindent and each of these two equivalent properties implies ${\cal G}_X = H_A^{n-1,\, n-1}(X,\,\R)$.

\end{Prop}

\noindent {\it Proof.} $(a)$\, The equivalence $(i)\Leftrightarrow(ii)$ follows by the standard duality and Hahn-Banach argument introduced in [Sul76] and used in various situations by several authors. Let $\Omega$ be a real $C^{\infty}$ form of bidegree $(n-1,\,n-1)$ on $X$. Then $\Omega$ is $d$-exact if and only if

\vspace{1ex}

\hspace{6ex}$\displaystyle\int\limits_X\Omega\wedge T = 0 \hspace{2ex} \mbox{for every real}\,\,d\mbox{-closed}\,\,(1,\,1)\mbox{-current}\,\,T\,\,\mbox{on}\,\,X,$

\vspace{1ex}

\noindent while $\Omega$ is positive definite if and only if

\vspace{1ex}

\hspace{6ex}$\displaystyle\int\limits_X\Omega\wedge T > 0 \hspace{2ex} \mbox{for every nonzero}\,\,(1,\,1)\mbox{-current}\,\,T\geq 0\,\,\mbox{on}\,\,X.$

\vspace{1ex}

\noindent It is thus clear that a form $\Omega$ as in $(i)$ and a current $T$ as in $(ii)$ cannot simultaneously exist. Thus $(i)\Rightarrow(ii)$. Conversely, if there is no $T$ as in $(ii)$, the set ${\cal E}$ of real $d$-closed $(1,\,1)$-currents $T$ on $X$ is disjoint from the set ${\cal C}$ of $(1,\,1)$-currents $T\geq 0$ on $X$ such that $\int_XT\wedge\gamma^{n-1}=1$ (where we have fixed an arbitrary smooth $(1,\,1)$-form $\gamma>0$ on $X$). Since ${\cal E}$ is a closed, convex subset of the locally convex space ${\cal D}_{\R}'$ of real $(1,\,1)$-currents on $X$, while ${\cal C}$ is a compact, convex subset of ${\cal D}_{\R}'$, by the Hahn-Banach separation theorem for locally convex spaces there must exist a linear functional on ${\cal D}_{\R}'$ that vanishes identically on ${\cal E}$ and is positive on ${\cal C}$ if ${\cal E}\cap{\cal C} =\emptyset$. This amounts to the existence of $\Omega$ as in $(i)$. The implication $(ii)\Rightarrow(i)$ is proved.

 We will now prove the equivalence ``$\mbox{not}\,(ii)\Leftrightarrow\mbox{not}\,(iii)$''. 

 Suppose there exists a non-trivial closed positive $(1,\, 1)$-current $T$ on $X$. If ${\cal G}_X$ degenerates, it contains the zero Aeppli $(n-1,\,n-1)$-class, so there exists a $C^{\infty}$ $(1,\,1)$-form $\omega>0$ on $X$ such that $\omega^{n-1}=\partial u + \bar\partial v$ for some forms $u,v$ of types $(n-2,\,n-1)$, resp. $(n-1,\,n-2)$. Thus, on the one hand, $\int_XT\wedge\omega^{n-1}>0$, while on the other hand Stokes's theorem would imply

$$\int\limits_XT\wedge\omega^{n-1} = \int\limits_XT\wedge(\partial u + \bar\partial v) = - \int\limits_X\partial T\wedge u - \int\limits_X\bar\partial T\wedge v = 0$$

\noindent since $\partial T=0$ and $\bar\partial T=0$ by the closedness assumption on $T$. This is a contradiction, so ${\cal G}_X$ cannot degenerate. We have thus proved the implication ``$\mbox{not}\,(ii)\Rightarrow\mbox{not}\,(iii)$''.
 
 Conversely, suppose that ${\cal G}_X\subsetneq H_A^{n-1,\,n-1}(X,\,\R)$. If no non-trivial closed positive $(1,\,1)$-current existed on $X$, then by the implication $(ii)\Rightarrow(i)$  proved above, there would exist a $d$-exact $C^{\infty}$ $(n-1,\,n-1)$-form $\Omega>0$ on $X$. Taking the $(n-1)^{st}$ root, there would exist a $C^{\infty}$ $(1,\,1)$-form $\omega>0$ on $X$ such that $\omega^{n-1}=\Omega$. Then $\omega^{n-1}\in\mbox{Im}\,d\subset\mbox{Im}\,\partial + \mbox{Im}\,\bar\partial$, hence $[\omega^{n-1}]_A=0$. However, $\omega$ is a Gauduchon (even a balanced) metric, so $[\omega^{n-1}]_A\in{\cal G}_X$. We would thus have $0\in{\cal G}_X$, hence ${\cal G}_X = H_A^{n-1,\,n-1}(X,\,\R)$, contradicting the assumption. This completes the proof of the implication ``$\mbox{not}\,(iii)\Rightarrow\mbox{not}\,(ii)$''.  

 The last statement in $(a)$ can be proved by contradiction. If $X$ were of class ${\cal C}$, then by the easy implication in Theorem 3.4 of [DP04] there would exist a K\"ahler current $T$ on $X$. However, any K\"ahler current is, in particular, a nonzero $d$-closed positive $(1,\, 1)$-current whose existence would violate $(ii)$. 

\vspace{1ex}

 To prove $(b)$, let us suppose that $H^2(X,\,\C)=0$. Then $H^{2n-2}(X,\,\C)=0$ by Poincar\'e duality, so for every balanced metric (if any) $\omega$ on $X$, $\omega^{n-1}$ must be $d$-exact, hence it must define a degenerate balanced structure on $X$. Thus, thanks to part $(a)$, $X$ is balanced if and only if there exists no nonzero $d$-closed $(1,\,1)$-current $T\geq 0$ on $X$. On the other hand, it was shown in [Pop09] that an arbitrary $X$ is sG if and only if there exists no nonzero $d$-exact $(1,\,1)$-current $T\geq 0$ on $X$. However, the assumption $H^2(X,\,\C)=0$ ensures that any $d$-closed current of degree $2$ is $d$-exact, so in this case the balanced and sG conditions on $X$ are characterised by the same property. This proves the equivalence in $(b)$.  

 The implication in $(b)$ follows from the above discussion\!\!: the assumption $H^2(X,\,\C)=0$ ensures that any balanced structure on $X$ is degenerate, while the existence of a degenerate balanced structure implies that the Gauduchon cone contains the zero Aeppli class, hence it must be the whole space $H_A^{n-1,\,n-1}(X,\,\R)$.  \hfill $\Box$

\vspace{2ex}

  We notice that the Gauduchon cone ${\cal G}_X$ and the sG cone ${\cal SG}_X$ cannot be simultaneously trivial, i.e. the implication holds\!\!:

\vspace{1ex}

\hspace{15ex} ${\cal G}_X=H^{n-1,\,n-1}_A(X,\,\R)\implies{\cal SG}_X\neq\emptyset.$

\vspace{1ex}

\noindent Indeed, if ${\cal G}_X=H^{n-1,\,n-1}_A(X,\,\R)$, then ${\cal SG}_X = \ker T\cap H_A^{n-1,\,n-1}(X,\,\R)$ is an $\R$ vector subspace of $H^{n-1,\,n-1}_A(X,\,\R)$, hence it contains at least the origin.

\vspace{2ex}

An immediate consequence of this and of Proposition \ref{Prop:exact-bal} is the following.

\begin{Cor}\label{Cor:Gcone-deg-conseq} If the Gauduchon cone ${\cal G}_X$ of a compact complex manifold $X$ degenerates, then $X$ is a strongly Gauduchon manifold but is not of class ${\cal C}$. 

\end{Cor}

 Recalling the implications ``$X$ is a class ${\cal C}$ manifold $\implies$ $X$ is a $\partial\bar\partial$-manifold $\implies$ $X$ is a strongly Gauduchon manifold'', the above corollary prompts the following question.

\begin{Question}\label{Question:Gcone-deg-ddbar} Do there exist $\partial\bar\partial$-manifolds $X$ whose Gauduchon cone ${\cal G}_X$ degenerates?

\end{Question}

 We notice that if such a manifold $X$ existed, it could not carry any pluriclosed metric. Indeed, it would have to carry a smooth $d$-exact $(n-1,\,n-1)$-form $\Omega>0$ by Proposition \ref{Prop:exact-bal} and $\Omega$ would have to be $\partial\bar\partial$-exact by the $\partial\bar\partial$-lemma. If a pluriclosed metric $\omega>0$ existed on $X$, then $\int_X\Omega\wedge\omega$ would have to be both positive and zero, a contradiction.

 A partial answer to Question \ref{Question:Gcone-deg-ddbar} may be contained in the discussion following Corollary 8.8 in [Fri91], although this is not clear to us since the notion of ``cohomologically K\"ahler'' manifold used there is said to be equivalent to that of manifold whose Fr\"olicher spectral sequence degenerates at $E_1$. If so, this notion is strictly weaker than our notion of a $\partial\bar\partial$-manifold. It would be very interesting to know whether the complex structure constructed in [Fri91] and [LT96] on $\sharp_k(S^3\times S^3)$ (for $k\geq 2$) satisfies the $\partial\bar\partial$ condition in the strong sense of the present work.

\vspace{2ex}  

 The duality (\ref{eqn:duality-BC-A-def}) between the Bott-Chern and Aeppli cohomologies can be restricted to various cones of cohomology classes. For example, if we consider the Bott-Chern K\"ahler cone of $X$, i.e. the open convex cone ${\cal K}_X\subset H^{1,\,1}_{BC}(X,\,\R)$ of Bott-Chern classes of K\"ahler metrics, we obviously have the following.

\begin{Obs}\label{Obs:K-cone_G-cone} The non-degenerate duality $H^{1,\,1}_{BC}(X,\,\C)\times H^{n-1,\,n-1}_A(X,\,\C) \rightarrow \C$ restricts to a {\bf positive} bilinear map

\vspace{1ex}

${\cal K}_X\times{\cal G}_X \longrightarrow \R, \hspace{3ex}  ([\omega]_{BC},\, [\gamma^{n-1}]_A) \longmapsto \int\limits_X\omega\wedge\gamma^{n-1}>0.$

\vspace{1ex}

\noindent In particular, ${\cal G}_X\subset({\cal K}_X)^{\mathrm{v}}$ and ${\cal K}_X\subset({\cal G}_X)^{\mathrm{v}}$, where for an open convex cone ${\cal C}$ in a finite-dimensional vector space $E$ we denote by ${\cal C}^{\mathrm{v}}$ the dual cone, i.e. the set of linear maps in $E^{\star}$ evaluating positively on every element in ${\cal C}$.

\end{Obs}

 It would be interesting to have an explicit description of the cone $({\cal G}_X)^{\mathrm{v}}\subset H^{1,\,1}_{BC}(X,\,\R)$ dual to the Gauduchon cone. The cone $({\cal G}_X)^{\mathrm{v}}$, which contains the Bott-Chern K\"ahler cone, is of course empty if ${\cal G}_X=H_A^{n-1,\,n-1}(X,\,\R)$, but when non-empty it may prove an efficient substitute for the K\"ahler cone when the latter is empty, so one may wonder if and to what extent it shares properties with it. 

 Similarly, recall Demailly's following definitions (cf. [Dem92]) of two other cones of Bott-Chern $(1,\,1)$-classes. The nef cone of $X$ is

\vspace{1ex}

 ${\cal NEF}_X:=\bigg\{\beta\in H^{1,\,1}_{BC}(X,\,\R)\,\slash\,\forall\,\varepsilon>0\,\,\exists\,\beta_{\varepsilon}\in\beta\,\,\mbox{smooth s.t.}\,\,\beta_{\varepsilon}\geq -\varepsilon\,\omega\bigg\},$

\vspace{1ex}

\noindent with $\omega>0$ a fixed $C^{\infty}$ $(1,\,1)$-form on $X$. If $X$ is K\"ahler, ${\cal NEF}_X$ is easily seen to be the closure of ${\cal K}_X$ (cf. [Dem92]). The pseudo-effective cone of $X$ is

\vspace{1ex}

\hspace{3ex} ${\cal E}_X:=\bigg\{[T]_{BC}\in H^{1,\,1}_{BC}(X,\,\R)\,\slash\, T\geq 0\,\,d\mbox{-closed}\,\,(1,\,1)\mbox{-current}\bigg\}.$

\vspace{1ex}

\noindent Clearly, ${\cal NEF}_X$ and ${\cal E}_X$ are {\it closed} convex cones (cf. [Dem92]) and ${\cal K}_X\subset{\cal NEF}_X\subset{\cal E}_X\subset H^{1,\,1}_{BC}(X,\,\R)$.

 Bearing in mind the duality between $H^{1,\,1}_{BC}(X,\,\R)$ and $H^{n-1,\,n-1}_A(X,\,\R)$, it seems natural to pursue in bidegree $(n-1,\,n-1)$ the analogy with the K\"ahler, nef and pseudo-effective cones of bidegree $(1,\,1)$. If the finite-dimensional vector space $H^{n-1,\,n-1}_A(X,\,\R)$ is endowed with the unique norm-induced topology, the closure in $H^{n-1,\,n-1}_A(X,\,\R)$ of the Gauduchon cone is the following closed convex cone

\vspace{1ex}

$\overline{\cal G}_X:=\bigg\{\alpha\in H^{n-1,\,n-1}_A(X,\,\R)\,\slash\,\forall\,\varepsilon>0\,\,\exists\,\Omega_{\varepsilon}\in\alpha\,\,\mbox{smooth s.t.}\,\,\Omega_{\varepsilon}\geq -\varepsilon\,\Omega\bigg\},$

\vspace{1ex}

\noindent where $\Omega>0$ is a fixed $C^{\infty}$ $(n-1,\,n-1)$-form on $X$ such that $\partial\bar\partial\Omega=0$. This follows immediately from the fact that a class $\alpha\in H^{n-1,\,n-1}_A(X,\,\R)$ is in the closure of ${\cal G}_X$ iff for every $\varepsilon>0$, $\alpha + \varepsilon\,[\Omega]_A\in{\cal G}_X$ (supposing that we have chosen $[\Omega]_A\neq 0\in H^{n-1,\,n-1}_A(X,\,\R)$\!; if $H^{n-1,\,n-1}_A(X,\,\R)=0$, everything is trivial). Clearly, by compactness of $X$, the definition of $\overline{\cal G}_X$ does not depend on the choice of $\Omega$. We can also define the cone ${\cal N}_X\subset H^{n-1,\,n-1}_A(X,\,\R)$\!:

\vspace{1ex}

\noindent ${\cal N}_X:=\bigg\{[U]_A\in H^{n-1,\,n-1}_A(X,\,\R)\,\slash\, U\geq 0\,\,\partial\bar\partial\mbox{-closed}\,\,(n-1,\,n-1)\mbox{-current}\bigg\}.$

\vspace{1ex}

\noindent It is clear that ${\cal G}_X\subset{\cal N}_X$, hence $\overline{\cal G}_X\subset\overline{\cal N}_X$. If ${\cal N}_X$ happens not to be closed (cf. Proposition \ref{Prop:closed-cones} below), we can replace it with its closure $\overline{\cal N}_X$. Thus we have cones ${\cal G}_X\subset\overline{\cal G}_X\subset\overline{\cal N}_X\subset H^{n-1,\,n-1}_A(X,\,\R)$. Meanwhile, if the Gauduchon cone degenerates (cf. Proposition \ref{Prop:exact-bal}), then ${\cal G}_X = \overline{\cal G}_X = {\cal N}_X = H^{n-1,\,n-1}_A(X,\,\R)$.

\begin{Prop}\label{Prop:closed-cones} Let $X$ be a compact complex manifold, $\mbox{dim}_{\C}=n$.

\noindent $(i)$ If $X$ is K\"ahler, the cone ${\cal N}_X$ is closed in $H^{n-1,\,n-1}_A(X,\,\R)$.

\noindent $(ii)$ If $X$ is of class ${\cal C}$, the inclusion $\overline{\cal G}_X\subset{\cal N}_X$ holds.

\end{Prop}

\noindent {\it Proof.} $(i)$\, Suppose that $X$ admits a K\"ahler metric $\omega$. If $(U_j)_{j\in\N}$ are $\partial\bar\partial$-closed positive $(n-1,\,n-1)$-currents such that the Aeppli classes $[U_j]_A$ converge to some class $\alpha\in H^{n-1,\,n-1}_A(X,\,\R)$ as $j\rightarrow\infty$, then $\int_XU_j\wedge\omega$ (depending only on $[U_j]_A$ thanks to $\omega$ being K\"ahler) converges to $\int_X\alpha\wedge\omega$, hence the positive currents $U_j$ are uniformly bounded in mass. Therefore, there exists a subsequence $U_{j_k}$ converging weakly to some $(n-1,\,n-1)$-current $U$. Then $U\geq 0$, $\partial\bar\partial U=0$ and $[U]_A=\alpha$, proving that $\alpha\in{\cal N}_X$. Thus ${\cal N}_X$ is closed.

\vspace{1ex}

\noindent $(ii)$\, Suppose that $X$ is of class ${\cal C}$. By [DP04], this amounts to the existence of a K\"ahler current $T$, i.e. a $d$-closed $(1,\,1)$-current such that $T\geq\delta\,\omega$ for some constant $\delta>0$ and some Hermitian metric $\omega>0$. Let $\alpha\in\overline{\cal G}_X$ and let $(\Omega_{\varepsilon})_{\varepsilon>0}$ be a family of $C^{\infty}$ $(n-1,\,n-1)$-forms in $\alpha$ such that $\Omega_{\varepsilon}\geq -\varepsilon\,\Omega$ for all $\varepsilon>0$ small. Then $\Omega_{\varepsilon} + \varepsilon\,\Omega\geq 0$ and $\int_X(\Omega_{\varepsilon} + \varepsilon\,\Omega)\wedge T = \int_X\Omega_{\varepsilon}\wedge T + \varepsilon\,\int_X\Omega\wedge T$ is bounded when $\varepsilon\downarrow 0$ since $\int_X\Omega_{\varepsilon}\wedge T$ is independent of $\varepsilon$ thanks to $[\Omega_{\varepsilon}]_A$ being independent of $\varepsilon$ and to $\partial T=0$ and $\bar\partial T=0$. Moreover, $\int_X(\Omega_{\varepsilon} + \varepsilon\,\Omega)\wedge T\geq \delta\,\int_X(\Omega_{\varepsilon} + \varepsilon\,\Omega)\wedge\omega\geq 0$, hence $\int_X(\Omega_{\varepsilon} + \varepsilon\,\Omega)\wedge\omega$ is bounded as $\varepsilon\downarrow 0$. Therefore the family $(\Omega_{\varepsilon} + \varepsilon\,\Omega)_{\varepsilon>0}$ admits a subsequence converging weakly to an $(n-1,\,n-1)$-current $U$ as $\varepsilon\downarrow 0$. We must have $U\geq 0$, $\partial\bar\partial U=0$ and $[U]_A=\alpha$, proving that $\alpha\in{\cal N}_X$. Thus $\overline{\cal G}_X\subset{\cal N}_X$ if $X$ is of class ${\cal C}$.   \hfill $\Box$

\vspace{2ex}

 It is natural to ask whether the K\"ahler assumption in $(i)$ or the class ${\cal C}$ assumption in $(ii)$ above may be relaxed. If we only suppose that ${\cal G}_X\subsetneq H^{n-1,\,n-1}_A(X,\,\R)$, Proposition \ref{Prop:exact-bal} ensures the existence of a nonzero $d$-closed $(1,\,1)$-current $T\geq 0$ for which the expressions $\int_X(\Omega_{\varepsilon} + \varepsilon\,\Omega)\wedge T$ in the proof of $(ii)$ in Proposition \ref{Prop:closed-cones} are still bounded when $\varepsilon\downarrow 0$. However, this is not enough to infer the existence of a weakly convergent subsequence of $(\Omega_{\varepsilon})_{\varepsilon>0}$. One may wonder what could be said if ``many'' $d$-closed positive $(1,\,1)$-currents $T$ existed on $X$. For example, if the algebraic dimension of $X$ is maximal (i.e. $a(X)=n$), then there are ``many'' divisors $D$ on $X$ inducing $d$-closed positive $(1,\,1)$-currents of integration $T=[D]$. However, $a(X)=n$ means that $X$ is Moishezon, hence $X$ is also of class ${\cal C}$ and we are in the situation of $(ii)$. 

 We now sum up the natural questions arising from the above considerations that we will hopefully take up in future work.

\begin{Question}\label{Question:overall} $(i)$\, Are the cones ${\cal NEF}_X$ and $\overline{\cal N}_X$, as well as the cones ${\cal E}_X$ and $\overline{\cal G}_X$, dual under the duality $H^{1,\,1}_{BC}(X,\,\C)\times H^{n-1,\,n-1}_A(X,\,\C) \rightarrow \C$\!?

\vspace{1ex}

 {\rm It is clear that we have inclusions ${\cal NEF}_X\subset\overline{\cal N}_X^{\mathrm{v}}$ and ${\cal E}_X\subset\overline{\cal G}_X^{\mathrm{v}}$, where for a closed convex cone ${\cal C}$ in a finite-dimensional vector space $E$ we denote by ${\cal C}^{\mathrm{v}}$ the dual cone, i.e. the set of linear maps in $E^{\star}$ evaluating non-negatively on every element in ${\cal C}$. It is also clear that if $X$ satisfies any of the equivalent conditions $(i), (ii), (iii)$ of part $(a)$ of Proposition \ref{Prop:exact-bal}, then ${\cal E}_X = \overline{\cal G}_X^{\mathrm{v}} = \{0\}$.} 

\vspace{1ex}

\noindent $(ii)$\, Can we define a notion of existence of ``many'' $d$-closed positive $(1,\,1)$-currents $T$ on $X$?

\vspace{1ex}

{\rm This might mean that the pseudo-effective cone ${\cal E}_X$ is ``maximal'' in some sense that has yet to be defined and also that the cone $\overline{\cal G}_X$ is ``minimal'' if these two cones are dual to each other. Any notion of ``minimality'' of ${\cal G}_X$ should be a strengthening of the property ${\cal SG}_X = {\cal G}_X$ which is necessary but not sufficient to ensure that $X$ is of class ${\cal C}$ (cf. Observation \ref{Obs:cones-equality}).}

\vspace{1ex}

\noindent $(iii)$\, If the answer to $(ii)$ is affirmative, does the following equivalence hold\!\!:

\vspace{1ex}

\noindent $X$ is of {\it class} ${\cal C} \iff$ there exist ``many'' $d$-closed positive $(1,\,1)$-currents $T$ on $X$\!\!?

\vspace{1ex}

 {\rm This would be the transcendental analogue of the standard characterisation of Moishezon manifolds as the compact complex manifolds carrying ``many'' divisors (i.e. having maximal algebraic dimension). 

If the answers to these questions turn out to be affirmative, then the {\it class} ${\cal C}$ manifolds will be precisely those compact complex manifolds whose Gauduchon cone is ``minimal''. If this proves to be the case, then the standard conjecture predicting that any deformation limit of {\it class} ${\cal C}$ manifolds is again of {\it class} ${\cal C}$ would follow since it will be seen below that the Gauduchon cone can only shrink or remain constant in the deformation limit.}

\end{Question}

 We shall now show that the Gauduchon cone behaves lower semicontinuously under holomorphic deformations of a $\partial\bar\partial$ complex structure. Let $\pi\,:\,{\cal X}\longrightarrow\Delta$ be a proper holomorphic submersion between complex manifolds. The question being local, we can assume that $\Delta\subset\C^m$ is an open ball containing the origin for some $m\in\N^{\star}$. All the fibres $X_t:=\pi^{-1}(t)$, $t\in\Delta$, are compact complex manifolds of equal dimensions $n$ and are $C^{\infty}$ diffeomorphic to a fixed $C^{\infty}$ manifold $X$, while the family of complex structures $(J_t)_{t\in\Delta}$ varies holomorphically with $t\in\Delta$. If we assume that $X_0$ is a $\partial\bar\partial$-manifold, the main result in [Wu06] ensures that $X_t$ is again a $\partial\bar\partial$-manifold for all $t\in\Delta$ sufficiently close to $0$. After possibly shrinking $\Delta$ about $0$, we may assume that this is the case for all $t\in\Delta$. Thus, by Theorem \ref{The:ddbar-Aeppli}, we have a Hodge-Aeppli decomposition on each fibre $X_t$ which in the case of the De Rham cohomology group $H^{2n-2}(X,\,\C)$ (necessarily independent of $t\in\Delta$) reads

$$H^{2n-2}(X,\,\C)\simeq H_A^{n,\,n-2}(X_t,\,\C)\oplus H_A^{n-1,\,n-1}(X_t,\,\C) \oplus H_A^{n-2,\,n}(X_t,\,\C), \hspace{3ex} t\in\Delta.$$
 
\noindent The $\partial\bar\partial$ assumption on the fibres $X_t$ ensures that the dimension of each of the spaces $H_A^{n,\,n-2}(X_t,\,\C)$, $H_A^{n-1,\,n-1}(X_t,\,\C)$ and $H_A^{n-2,\,n}(X_t,\,\C)$ is independent of $t\in\Delta$. Therefore the ellipticity of the Aeppli Laplacians $\Delta_A^{(t)}$ (defined by any smooth family of Hermitian metrics $(\omega_t)_{t\in\Delta}$ on the fibres $(X_t)_{t\in\Delta}$) and the Kodaira-Spencer theory [KS60] imply that

$$\Delta\ni t\mapsto H_A^{n-1,\,n-1}(X_t,\,\C)$$

\noindent and its analogues in bidegrees $(n,\,n-2)$, $(n-2,\,n)$ are $C^{\infty}$ vector bundles, while the projections of $H^{2n-2}(X,\,\C)$ on
$H_A^{n,\,n-2}(X_t,\,\C)$, $H_A^{n-1,\,n-1}(X_t,\,\C)$ and $H_A^{n-2,\,n}(X_t,\,\C)$ vary in a $C^{\infty}$ way with $t\in\Delta$. Thus, composing the canonical injection $H_A^{n-1,\,n-1}(X_0,\,\C)\hookrightarrow H^{2n-2}(X,\,\C)$ of Theorem \ref{The:ddbar-Aeppli} with the canonical projection $H^{2n-2}(X,\,\C)\twoheadrightarrow H_A^{n-1,\,n-1}(X_t,\,\C)$ induced by the Hodge-Aeppli decomposition, we get a linear map

\begin{equation}\label{eqn:A_tdef}A_t\,:\,H_A^{n-1,\,n-1}(X_0,\,\C)\longrightarrow H_A^{n-1,\,n-1}(X_t,\,\C), \hspace{2ex} t\in\Delta,\end{equation}

\noindent that depends in a $C^{\infty}$ way on $t$. Since $A_0$ is the identity of $H_A^{n-1,\,n-1}(X_0,\,\C)$, $A_t$ must be an isomorphism of complex vector spaces for all $t\in\Delta$ after possibly further shrinking $\Delta$ about $0$. 

 The isomorphisms $A_t$ in (\ref{eqn:A_tdef}) can be used to compare ${\cal G}_{X_0}$ with ${\cal G}_{X_t}$.

\begin{The}\label{The:Gcone-semicont} The Gauduchon cones $({\cal G}_{X_t})_{t\in\Delta}$ of the fibres $(X_t)_{t\in\Delta}$ of any holomorphic family of $\partial\bar\partial$-manifolds satisfy the following semi-continuity property. For every $[\omega_0^{n-1}]_A\in{\cal G}_{X_0}$, there exists $\varepsilon>0$ (depending on $[\omega_0^{n-1}]_A$) such that

$$A_t([\omega_0^{n-1}]_A)\in{\cal G}_{X_t}  \hspace{2ex} \mbox{for all}\hspace{1ex} t\in\Delta\hspace{1ex}\mbox{with}\hspace{1ex} |t|<\varepsilon.$$

\end{The}

 In other words, if we identify every ${\cal G}_{X_t}$ with its image in $H^{2n-2}(X,\,\C)$ under the canonical injection $H_A^{n-1,\,n-1}(X_t,\,\C)\hookrightarrow H^{2n-2}(X,\,\C)$ for every $t$, the Gauduchon cone of $X_0$ is contained in the limit as $t$ approaches $0$ of the Gauduchon cones of $X_t$. So in a sense the Gauduchon cone can only shrink or remain constant on the limit fibre. Note that if we do not make the $\partial\bar\partial$ assumption on the fibres $(X_t)_{t\in\Delta}$, the picture may change\!: we may have $\mbox{dim}\,H_A^{n-1,\,n-1}(X_0,\,\C)>\mbox{dim}\,H_A^{n-1,\,n-1}(X_t,\,\C)$ for $t\neq 0$, so in this case the dimension of ${\cal G}_{X_0}$ as a complex manifold $(=$ the dimension of $H_A^{n-1,\,n-1}(X_0,\,\C)$ as a vector space since ${\cal G}_{X_0}\subset H_A^{n-1,\,n-1}(X_0,\,\C)$ is an open subset$)$ is strictly larger than the dimension of ${\cal G}_{X_t}$ as a complex manifold for $t\neq 0$. 

 Before proving Theorem \ref{The:Gcone-semicont}, we notice the following.

\begin{Lem}\label{Lem:inj-proj} Let $X$ be a compact complex manifold, $\mbox{dim}_{\C}X=n$. Fix an arbitrary smooth $(2n-2)$-form $\Omega$ on $X$ such that $d\Omega=0$.

\noindent $(i)$\, If $\Omega=\Omega^{n,\,n-2} + \Omega^{n-1,\,n-1} + \Omega^{n-2,\,n}$ is the splitting into components of pure types, then

$$\partial\bar\partial\Omega^{n,\,n-2}=0, \hspace{2ex} \partial\bar\partial\Omega^{n-1,\,n-1}=0, \hspace{2ex} \partial\bar\partial\Omega^{n-2,\,n}=0.$$

\noindent $(ii)$\, Suppose that $X$ is a $\partial\bar\partial$-manifold. Then

$$\{\Omega\} = [\Omega^{n,\,n-2}]_A + [\Omega^{n-1,\,n-1}]_A + [\Omega^{n-2,\,n}]_A,$$

\noindent where $\{\Omega\}\in H^{2n-2}(X,\,\C)$ denotes the De Rham class of $\Omega$, while $[\Omega^{p,\,q}]_A$ denotes the image in $H^{2n-2}(X,\,\C)$ of the Aeppli class of $\Omega^{p,\,q}$ under the canonical injection $H_A^{p,\,q}(X,\,\C)\hookrightarrow H^{2n-2}(X,\,\C)$ defined by the $\partial\bar\partial$ assumption on $X$ (cf. Theorem \ref{The:ddbar-Aeppli}) for all $(p,\,q)\in\{(n,\,n-2), (n-1,\,n-1), (n-2,\,n)\}$. (Thus we denote by the same symbol an Aeppli class and its canonical image into De Rham cohomology.)

\end{Lem}

\noindent {\it Proof.} The form $d\Omega$ is of degree $(2n-1)$, so it has two pure-type components of bidegrees $(n,\,n-1)$, resp. $(n-1,\,n)$. Thus $d\Omega=0$ amounts to the vanishing of each of these\!:

\begin{equation}\label{eqn:nn-1-vanishing}(a)\,\partial\Omega^{n-1,\,n-1} + \bar\partial\Omega^{n,\,n-2}=0  \hspace{2ex} \mbox{and} \hspace{2ex} (b)\, \partial\Omega^{n-2,\,n} + \bar\partial\Omega^{n-1,\,n-1}=0.\end{equation}

\noindent Applying $\partial$ in $(b)$ (or $\bar\partial$ in $(a)$), we get $\partial\bar\partial\Omega^{n-1,\,n-1}=0$. Now, $\Omega^{n,\,n-2}$ is $\partial$-closed and $\Omega^{n-2,\,n}$ is $\bar\partial$-closed for obvious bidegree reasons, hence they are also $\partial\bar\partial$-closed. This proves $(i)$.

 To prove $(ii)$, we have to spell out the canonical images of the Aeppli classes $[\Omega^{p,\,q}]_A$ into De Rham cohomology as in the proof of Theorem \ref{The:ddbar-Aeppli}.

 In the case of $\Omega^{n-1,\,n-1}$, we need forms $\xi$, $\eta$ of bidegrees $(n-2,\,n-1)$, resp. $(n-1,\,n-2)$, such that $d(\partial\xi + \Omega^{n-1,\,n-1} + \bar\partial\eta)=0$, which amounts to

$$\partial\bar\partial\xi = \bar\partial\Omega^{n-1,\,n-1} \hspace{2ex} \mbox{and} \hspace{2ex} \partial\bar\partial\eta = -\partial\Omega^{n-1,\,n-1}.$$

\noindent If we fix a Hermitian metric $\omega$ on $X$ and choose $\xi$ and $\eta$ to be the solutions of minimal $L^2$ norms of these $\partial\bar\partial$ equations, formula (\ref{eqn:ddbar-eq-solformula1}) of Theorem \ref{The:min-sol-ddbar} gives

$$\xi = (\partial\bar\partial)^{\star}\Delta_{BC}^{-1}\bigg(\bar\partial\Omega^{n-1,\,n-1}\bigg) \hspace{2ex} \mbox{and} \hspace{2ex} \eta = - (\partial\bar\partial)^{\star}\Delta_{BC}^{-1}\bigg(\partial\Omega^{n-1,\,n-1}\bigg).$$

\noindent The form $\Gamma^{n-1,\,n-1}:=\partial\xi + \Omega^{n-1,\,n-1} + \bar\partial\eta$ constructed in this way reads

$$\Gamma^{n-1,\,n-1}:=\partial(\partial\bar\partial)^{\star}\Delta_{BC}^{-1}\bigg(\bar\partial\Omega^{n-1,\,n-1}\bigg) + \Omega^{n-1,\,n-1} - \bar\partial(\partial\bar\partial)^{\star}\Delta_{BC}^{-1}\bigg(\partial\Omega^{n-1,\,n-1}\bigg),$$

\noindent is of bidegree $(n-1,\,n-1)$, $d$-closed and Aeppli cohomologous to $\Omega^{n-1,\,n-1}$. Thus the canonical image of $[\Omega^{n-1,\,n-1}]_A = [\Gamma^{n-1,\,n-1}]_A\in H_A^{n-1,\,n-1}(X,\,\C)$ into $H^{2n-2}(X,\,\C)$ is the De Rham class $\{\Gamma^{n-1,\,n-1}\}$.

 Running the same procedure for $\Omega^{n,\,n-2}$ and $\Omega^{n-2,\,n}$, we get $d$-closed forms

\begin{eqnarray}\nonumber\Gamma^{n,\,n-2} & = & \partial(\partial\bar\partial)^{\star}\Delta_{BC}^{-1}\bigg(\bar\partial\Omega^{n,\,n-2}\bigg) + \Omega^{n,\,n-2}, \hspace{2ex}\mbox{of bidegree}\hspace{1ex} (n,\,n-2),\\ 
\nonumber\Gamma^{n-2,\,n} & = & \Omega^{n-2,\,n} - \bar\partial(\partial\bar\partial)^{\star}\Delta_{BC}^{-1}\bigg(\partial\Omega^{n-2,\,n}\bigg), \hspace{2ex}\mbox{of bidegree}\hspace{1ex} (n-2,\,n),\end{eqnarray}

\noindent that are Aeppli cohomologous to $\Omega^{n,\,n-2}$, resp. $\Omega^{n-2,\,n}$. To finish the proof of $(ii)$, it remains to prove the following identity of De Rham classes

\begin{equation}\label{eqn:DeRhan-identity}\{\Omega\} = \{\Gamma^{n,\,n-2} + \Gamma^{n-1,\,n-1} + \Gamma^{n-2,\,n}\}.\end{equation}

\noindent We see that $\Gamma^{n,\,n-2} + \Gamma^{n-1,\,n-1} + \Gamma^{n-2,\,n} = \Omega + \partial\alpha + \bar\partial\beta$, where

$$\alpha = (\partial\bar\partial)^{\star}\Delta_{BC}^{-1}\bigg(\bar\partial\Omega^{n,\,n-2} + \bar\partial\Omega^{n-1,\,n-1}\bigg) \hspace{1ex} \mbox{and} \hspace{1ex} \beta = - (\partial\bar\partial)^{\star}\Delta_{BC}^{-1}\bigg(\partial\Omega^{n-2,\,n} + \partial\Omega^{n-1,\,n-1}\bigg).$$

\noindent Now, formulae (\ref{eqn:nn-1-vanishing}) show that $\alpha=\beta$. Indeed, $(a)$ and $(b)$ add up to $\bar\partial\Omega^{n,\,n-2} + \bar\partial\Omega^{n-1,\,n-1} = - (\partial\Omega^{n-2,\,n} + \partial\Omega^{n-1,\,n-1})$. We get $\partial\alpha + \bar\partial\beta = \partial\alpha + \bar\partial\alpha =d\alpha$, hence

$$\Gamma^{n,\,n-2} + \Gamma^{n-1,\,n-1} + \Gamma^{n-2,\,n} = \Omega + d\alpha.$$

\noindent This proves (\ref{eqn:DeRhan-identity}) and completes the proof of the lemma.  \hfill $\Box$

\vspace{3ex}

\noindent {\it Proof of Theorem \ref{The:Gcone-semicont}.} Let $[\omega_0^{n-1}]_A\in H_A^{n-1,\,n-1}(X_0,\,\C)$ be an arbitrary element in ${\cal G}_{X_0}$, where $\omega_0>0$ is a Gauduchon metric on $X_0$. 

Thanks to the $\partial\bar\partial$ assumption, we can find forms $u_0$ and $v_0$ of respective $J_0$-types $(n-2,\,n-1)$ and $(n-1,\,n-2)$ such that

$$\Omega:=\partial_0u_0 + \omega_0^{n-1} + \bar\partial_0v_0$$

\noindent is $d$-closed. Let $(\Omega_t^{n-1,\,n-1})_{t\in\Delta}$ be the $C^{\infty}$ family of components of $\Omega$ of $J_t$-types $(n-1,\,n-1)$. By $(i)$ of Lemma \ref{Lem:inj-proj}, we have $\partial_t\bar\partial_t\Omega_t^{n-1,\,n-1}=0$ for all $t$. We extend $u_0$ and $v_0$ in an arbitrary way to $C^{\infty}$ families $(u_t)_{t\in\Delta}$ and $(v_t)_{t\in\Delta}$ of forms of $J_t$-types $(n-2,\,n-1)$ and resp. $(n-1,\,n-2)$ and we set

$$\Lambda_t:=\Omega_t^{n-1,\,n-1} - \partial_tu_t - \bar\partial_tv_t, \hspace{2ex} t\in\Delta.$$ 

\noindent It is clear that $\partial_t\bar\partial_t\Lambda_t=0$, that $[\Lambda_t]_A = [\Omega_t^{n-1,\,n-1}]_A$ and that the family of forms $(\Lambda_t)_{t\in\Delta}$ of $J_t$-types $(n-1,\,n-1)$ depends in a $C^{\infty}$ way on $t\in\Delta$ and 

$$\Lambda_0 = \Omega - \partial_0u_0 - \bar\partial_0v_0 = \omega_0^{n-1}>0$$

\noindent since $\Omega$ is of type $(n-1,\,n-1)$ for $J_0$, so $\Omega_0^{n-1,\,n-1}=\Omega$. By the continuity of the family $(\Lambda_t)_{t\in\Delta}$, the strict positivity of $\Lambda_0$ implies the strict positivity of $\Lambda_t$ for all $t\in\Delta$ sufficiently close to $0$. Thus we can extract the $(n-1)^{st}$ root\!: for every $t$ close to $0$, there exists a unique positive definite smooth form $\omega_t$ of $J_t$-type $(1,\,1)$ such that $\omega_t^{n-1}=\Lambda_t$. Every such $\omega_t$ is thus a Gauduchon metric on $X_t$ and we have

$$A_t([\omega_0^{n-1}]_A)=A_t([\Omega]_A)\stackrel{(a)}{=}[\Omega_t^{n-1,\,n-1}]_A=[\Lambda_t]_A=[\omega_t^{n-1}]_A\in{\cal G}_{X_t},$$

\noindent where the identity $(a)$ above follows from $(ii)$ of Lemma \ref{Lem:inj-proj} applied to the $\partial\bar\partial$ complex structure $J_t$.  \hfill $\Box$

\vspace{2ex}

 A more precise description of the variation of the Gauduchon cone ${\cal G}_X$ under deformations of $X$ may be possible after singling out a special representative for every Aeppli-Gauduchon class by solving equation $(\star)$.

\section{Proof of uniqueness in equation $(\star)$}\label{section:uniqueness}

 We start by proving the uniqueness of solutions to equation $(\star\star)$ subject to (\ref{eqn:conditions-second}) when the given $\omega$ is an arbitrary Hermitian metric.

\begin{Prop}\label{Prop:uniqueness-second} Let $(X,\,\omega)$ be a compact Hermitian manifold, $\mbox{dim}_{\C}X=n\geq 2$. Suppose that for real-valued $C^{\infty}$ functions $\varphi_1$ and $\varphi_2$ on $X$ we have $\omega^{n-1} + i\partial\bar\partial\varphi_l\wedge\omega^{n-2} > 0$ (for $l=1,2$) and

\begin{equation}\label{eqn:det-equality}\bigg[\bigg( \omega^{n-1} + i\partial\bar\partial\varphi_1\wedge\omega^{n-2}\bigg)^{\frac{1}{n-1}}\bigg]^n = \bigg[\bigg( \omega^{n-1} + i\partial\bar\partial\varphi_2\wedge\omega^{n-2}\bigg)^{\frac{1}{n-1}}\bigg]^n.\end{equation}

\noindent Then the function $\varphi_1 - \varphi_2$ is constant on $X$.

\end{Prop}

\vspace{2ex}

  We begin on a few preliminary calculations that will prove useful later on. The symbol $\Lambda=\Lambda_{\omega}$ will stand for the formal adjoint of the Lefschetz operator $L_{\omega}=\omega\wedge\cdot$ of multiplication by the Hermitian metric $\omega$, while $\Delta_{\omega}\varphi:=\Lambda_{\omega}(i\partial\bar\partial\varphi)$ will denote the (non-positive) Laplacian associated with $\omega$ on real-valued $C^2$ functions $\varphi$ on $X$. On $(1,\,1)$-forms, $\Lambda_{\omega}$ coincides with the trace w.r.t. $\omega$ denoted by $\mbox{tr}_{\omega}$, so $\Lambda_{\omega}$ and $\mbox{tr}_{\omega}$ will be used interchangeably.

\begin{Lem}\label{Lem:11-prim} Let $(X,\,\omega)$ be a compact Hermitian manifold, $\mbox{dim}_{\C}X=n$.

\noindent $(i)$\, For any smooth $(1,\,1)$-form $\alpha$ on $X$, the Lefschetz decomposition of $\alpha$ w.r.t. $\omega$ (into forms of bidegree $(1,\,1)$) reads

\begin{equation}\label{eqn:Lefschetz-decomp}\alpha = \alpha_{prim} + \frac{1}{n}\,(\Lambda_{\omega}\alpha)\,\omega,\end{equation}

\noindent where the primitive part $\alpha_{prim}$ of $\alpha$ is defined by either of the equivalent conditions: $\Lambda_{\omega}\alpha_{prim}=0$ or $\alpha_{prim}\wedge\omega^{n-1}=0$.

\noindent $(ii)$\, In particular, if $\star=\star_{\omega}$ is the Hodge star operator defined by $\omega$, we have

\begin{equation}\label{eqn:star-alpha-n-2}\star\bigg(\alpha\wedge\frac{\omega^{n-2}}{(n-2)!}\bigg) = -\alpha + (\Lambda_{\omega}\alpha)\,\omega.\end{equation}

\noindent Hence, if $\alpha=i\partial\bar\partial\varphi$ for some real-valued function $\varphi$, then

\begin{equation}\label{eqn:star-iddbar_phi-n-2}\star\bigg(i\partial\bar\partial\varphi\wedge\frac{\omega^{n-2}}{(n-2)!}\bigg) = -i\partial\bar\partial\varphi + (\Delta_{\omega}\varphi)\,\omega.\end{equation}

\noindent $(iii)$\, Still denoting $\star=\star_{\omega}$, for any smooth $(n-1,\,n-1)$-form $\Gamma$, we have

\begin{equation}\label{eqn:trace-star}\mbox{tr}_{\omega_{n-1}}\Gamma = \mbox{tr}_{\omega}(\star\Gamma).\end{equation}

\noindent $(iv)$\, For any real-valued $C^2$ function $\varphi$ on $X$, we have

\begin{equation}\label{eqn:trace-iddbar_n-2}\mbox{tr}_{\omega_{n-1}}\bigg(i\partial\bar\partial\varphi\wedge\frac{\omega^{n-2}}{(n-2)!}\bigg) = (n-1)\,\Delta_{\omega}\varphi,\end{equation}

\noindent where $\mbox{tr}_{\omega_{n-1}}$ denotes the trace w.r.t. $\omega_{n-1}:=\frac{\omega^{n-1}}{(n-1)!}$ of the $(n-1,\,n-1)$-form to which it applies.

\end{Lem}

\noindent {\it Proof.} $(i)$\, By the Lefschetz decomposition, any $\alpha\in C^{\infty}_{1,\,1}(X,\,\C)$ splits as $\alpha = \alpha_{prim} + f\,\omega$ for a unique primitive $(1,\,1)$-form $\alpha_{prim}$ and a unique function $f$ on $X$. Applying $\Lambda_{\omega}$ and using $\Lambda_{\omega}\alpha_{prim}=0$, $\Lambda_{\omega}\omega=n$, we get $(\ref{eqn:Lefschetz-decomp})$.

 $(ii)$\, It is well known that for any primitive $(1,\,1)$-form $\alpha_{prim}$, we have

\begin{equation}\label{eqn:star-alpha-prim}\star\alpha_{prim} = -\alpha_{prim}\wedge\frac{\omega^{n-2}}{(n-2)!} \hspace{2ex} \mbox{i.e.} \hspace{2ex} \star\bigg(\alpha_{prim}\wedge\frac{\omega^{n-2}}{(n-2)!}\bigg) = - \alpha_{prim},\end{equation}

\noindent (i.e. the bidegree $(1,\,1)$ case of the formula $\star\, v = (-1)^{k(k+1)/2}\, i^{p-q}\, \frac{\omega^{n-p-q}\wedge v}{(n-p-q)!}$ for any {\it primitive} $(p,\,q)$-form $v$, where $k=p+q$ -- see e.g. [Voi02, Proposition 6.29, p. 150]). On the other hand, $\star(\omega^{n-1}/(n-1)!) = \omega$, so using (\ref{eqn:Lefschetz-decomp}) we get

$$\star\bigg(\alpha\wedge\frac{\omega^{n-2}}{(n-2)!}\bigg) = -\alpha_{prim} + \frac{n-1}{n}\,(\Lambda_{\omega}\alpha)\,\omega = -\alpha + (\Lambda_{\omega}\alpha)\,\omega.$$

$(iii)$\, To prove the pointwise identity $(\ref{eqn:trace-star})$, we fix an arbitrary point $x\in X$ and choose local holomorphic coordinates $z_1,\dots , z_n$ about $x$ such that 

$$\omega(x) = \sum\limits_{j=1}^nidz_j\wedge d\bar{z}_j \hspace{2ex} \mbox{and} \hspace{2ex} \Gamma(x) = \sum\limits_{j=1}^n\Gamma_j\,\widehat{idz_j\wedge d\bar{z}_j},$$

\noindent where for all $j=1,\dots , n$, we denote by $\widehat{idz_j\wedge d\bar{z}_j}$ the $(n-1,\,n-1)$-form $idz_1\wedge d\bar{z}_1\wedge\dots\wedge\widehat{(idz_j\wedge d\bar{z}_j)}\wedge\dots\wedge idz_n\wedge d\bar{z}_n$ (where\,\, $\widehat{ }$\,\, indicates a missing factor). It is clear that

\begin{equation}\label{eqn:star_n-1-coord}\star(\widehat{idz_j\wedge d\bar{z}_j}) = idz_j\wedge d\bar{z}_j \hspace{3ex} \mbox{at}\hspace{1ex} x.\end{equation}

\noindent Indeed, if we denote $dV_{\omega}=\omega^n/n!$, we have $(\widehat{idz_k\wedge d\bar{z}_k})\wedge(idz_j\wedge d\bar{z}_j) = \delta_{jk}\,idz_1\wedge d\bar{z}_1\wedge\dots\wedge idz_n\wedge d\bar{z}_n = \delta_{jk}\,dV_{\omega}(x) = \langle\widehat{idz_k\wedge d\bar{z}_k},\, \widehat{idz_j\wedge d\bar{z}_j}\rangle \,dV_{\omega}(x)$ for all $j,k=1, \dots , n$, where $\langle\,\,,\,\,\rangle$ stands for the pointwise scalar product defined by $\omega$ at $x$ and $\delta_{jk}$ is the Kronecker delta. 

Thus $(\star\Gamma)(x) = \sum\limits_{j=1}^n\Gamma_j\,idz_j\wedge d\bar{z}_j$. Meanwhile, $\omega_{n-1}(x) = \sum\limits_{j=1}^n \widehat{idz_j\wedge d\bar{z}_j}$, hence $(\mbox{tr}_{\omega_{n-1}}\Gamma)(x) = \sum\limits_{j=1}^n\Gamma_j = \mbox{tr}_{\omega}(\star\Gamma)(x)$, which proves (\ref{eqn:trace-star}). 

 $(iv)$\, We now use (\ref{eqn:trace-star}) and (\ref{eqn:star-iddbar_phi-n-2}) to get

$$\mbox{tr}_{\omega_{n-1}}\bigg(i\partial\bar\partial\varphi\wedge\frac{\omega^{n-2}}{(n-2)!}\bigg) = \mbox{tr}_{\omega}\star\bigg(i\partial\bar\partial\varphi\wedge\frac{\omega^{n-2}}{(n-2)!}\bigg) = -\Delta_{\omega}\varphi + n\,\Delta_{\omega}\varphi,$$

\noindent which proves $(\ref{eqn:trace-iddbar_n-2})$.    \hfill $\Box$

\vspace{3ex}

 If $\Gamma>0$ is an $(n-1,\,n-1)$-form for which the local coordinates have been chosen at a given point $x$ such that $\frac{\Gamma(x)}{(n-1)!} = \sum\limits_{j=1}^n\Gamma_j\,\widehat{idz_j\wedge d\bar{z}_j}$, then its $(n-1)^{st}$ root $\gamma=\Gamma^{\frac{1}{n-1}}$ is given at $x$ by

\begin{equation}\label{eqn:root-coeff}\gamma(x) = \sum\limits_{j=1}^n\gamma_j\,idz_j\wedge d\bar{z}_j, \hspace{2ex} \mbox{where} \hspace{2ex} \gamma_j = \frac{(\Gamma_1\dots \Gamma_n)^{\frac{1}{n-1}}}{\Gamma_j}, \hspace{2ex} j=1,\dots , n.\end{equation}

\noindent The $\gamma_j$'s are well-defined since $\Gamma_j>0$ for all $j$. In particular, we see that the determinants (which make intrinsic sense) are related by

\begin{equation}\label{eqn:determinants}\det\,(\Gamma^{\frac{1}{n-1}}) = \frac{1}{((n-1)!)^{\frac{n}{n-1}}}\,(\det\Gamma)^{\frac{1}{n-1}}.\end{equation}

 Now, let $\varphi_1$ and $\varphi_2$ be real-valued functions on $X$ as in the statement of Proposition \ref{Prop:uniqueness-second}. Fix an arbitrary point $x\in X$ and choose local holomorphic coordinates $z_1, \dots , z_n$ about $x$ such that

\vspace{1ex}

$\omega(x) = \sum\limits_{j=1}^nidz_j\wedge d\bar{z}_j \hspace{2ex} \mbox{and} \hspace{2ex} i\partial\bar\partial\varphi_l(x) = \sum\limits_{j=1}^n\lambda^{(l)}_j\,idz_j\wedge d\bar{z}_j, \hspace{2ex} l=1,2.$ 

\vspace{1ex}

\noindent Straightforward calculations give $\omega^{n-1}(x) = (n-1)!\,\sum\limits_{j=1}^n\widehat{idz_j\wedge d\bar{z}_j}$ and

\vspace{1ex}

$i\partial\bar\partial\varphi_l\wedge\omega^{n-2}(x) = (n-2)!\,\sum\limits_{j=1}^n(\lambda_1^{(l)} + \dots + \lambda_n^{(l)} - \lambda_j^{(l)})\,\widehat{idz_j\wedge d\bar{z}_j}, \hspace{2ex} l=1,2.$

\vspace{1ex}

\noindent Hence, if we set $\mu_j^{(l)}:=1 + \lambda_j^{(l)}$ and $\xi_j^{(l)} = \mu_1^{(l)} + \dots + \mu_n^{(l)} - \mu_j^{(l)}$ for $1\leq j\leq n$ and $l=1,2$, we get at $x$\!:

\begin{equation}\label{eqn:pre-root-forms}\omega^{n-1} + i\partial\bar\partial\varphi_l\wedge\omega^{n-2} = (n-2)!\,\sum\limits_{j=1}^n\xi_j^{(l)}\,\widehat{idz_j\wedge d\bar{z}_j}, \hspace{2ex} l=1,2.\end{equation}

\noindent Using (\ref{eqn:root-coeff}), we see that  at $x$ the roots for $l=1,2$ read

\vspace{1ex}

$\bigg(\frac{\omega^{n-1} + i\partial\bar\partial\varphi_l\wedge\omega^{n-2}}{(n-2)!}\bigg)^{\frac{1}{n-1}} = \frac{1}{((n-1)!)^{\frac{1}{n-1}}}\,\bigg(\prod\limits_{k=1}^n\xi_k^{(l)}\bigg)^{\frac{1}{n-1}}\,\sum\limits_{j=1}^n\frac{1}{\xi_j^{(l)}}\,idz_j\wedge d\bar{z}_j,$

\vspace{1ex}

\noindent while using (\ref{eqn:star_n-1-coord}), we have at $x$\!:

\begin{equation}\label{eqn:star-forms}\star(\omega^{n-1} + i\partial\bar\partial\varphi_l\wedge\omega^{n-2}) =  (n-2)!\,\sum\limits_{j=1}^n\xi_j^{(l)}\,idz_j\wedge d\bar{z}_j, \hspace{2ex} l=1,2.\end{equation}

\noindent {\it Proof of Proposition \ref{Prop:uniqueness-second}.} For real-valued functions $\varphi_1$ and $\varphi_2$ on $X$ as in the statement of Proposition \ref{Prop:uniqueness-second}, we consider the positive definite $(1,\,1)$-forms

$$\gamma_1:=(\omega^{n-1} + i\partial\bar\partial\varphi_1\wedge\omega^{n-2})^{\frac{1}{n-1}}  \hspace{2ex} \mbox{and} \hspace{2ex} \gamma_2:=(\omega^{n-1} + i\partial\bar\partial\varphi_2\wedge\omega^{n-2})^{\frac{1}{n-1}}.$$

\noindent Hypothesis (\ref{eqn:det-equality}) translates to the following sequence of equivalent identities

\begin{eqnarray}\label{gamma_1,2-eq1}\nonumber\gamma_1^n = \gamma_2^n & \iff & \det(\gamma_1) = \det(\gamma_2) \stackrel{(a)}{\iff} \det(\gamma_1^{n-1}) = \det(\gamma_2^{n-1})\\
\nonumber & \iff & \det(\omega^{n-1} + i\partial\bar\partial\varphi_1\wedge\omega^{n-2}) = \det(\omega^{n-1} + i\partial\bar\partial\varphi_2\wedge\omega^{n-2})\\
\nonumber & \stackrel{(b)}{\iff} & \det\bigg(\star(\omega^{n-1} + i\partial\bar\partial\varphi_1\wedge\omega^{n-2})\bigg) = \det\bigg(\star(\omega^{n-1} + i\partial\bar\partial\varphi_2\wedge\omega^{n-2})\bigg)\\
\nonumber & \iff & \bigg(\star(\omega^{n-1} + i\partial\bar\partial\varphi_1\wedge\omega^{n-2})\bigg)^n = \bigg(\star(\omega^{n-1} + i\partial\bar\partial\varphi_2\wedge\omega^{n-2})\bigg)^n\\
\nonumber & \iff & \star\bigg(i\partial\bar\partial(\varphi_1-\varphi_2)\wedge\omega^{n-2}\bigg)\wedge\rho^{n-1}=0\\
     & \stackrel{(c)}{\iff} & \bigg(-i\partial\bar\partial(\varphi_1-\varphi_2) + \Delta_{\omega}(\varphi_1-\varphi_2)\,\omega\bigg)\wedge\rho^{n-1}=0,\end{eqnarray}

\noindent where $(a)$ follows from (\ref{eqn:determinants}), $(b)$ follows from comparing (\ref{eqn:pre-root-forms}) and (\ref{eqn:star-forms}), $(c)$ follows from (\ref{eqn:star-iddbar_phi-n-2}), while $\rho$ denotes the smooth, positive definite $(1,\,1)$-form that is the $(n-1)^{st}$ root of the smooth, positive definite $(n-1,\,n-1)$-form

\vspace{1ex}

$\Omega:=\sum\limits_{p=1}^n\bigg(\star(\omega^{n-1} + i\partial\bar\partial\varphi_1\wedge\omega^{n-2})\bigg)^{n-p}\wedge\bigg(\star(\omega^{n-1} + i\partial\bar\partial\varphi_2\wedge\omega^{n-2})\bigg)^{p-1}.$

\vspace{1ex}

\noindent Further transforming (\ref{gamma_1,2-eq1}), we get

\begin{eqnarray}\label{gamma_1,2-eq2}\nonumber\gamma_1^n = \gamma_2^n & \iff & \Delta_{\omega}(\varphi_1-\varphi_2)\,\omega\wedge\rho^{n-1} - \Delta_{\rho}(\varphi_1-\varphi_2)\,\frac{\rho^n}{n}=0\\
\nonumber   & \iff & \bigg(\Delta_{\omega}(\varphi_1-\varphi_2)\,\omega - \frac{1}{n}\,\Delta_{\rho}(\varphi_1-\varphi_2)\,\rho\bigg)\wedge\rho^{n-1}=0\\
\nonumber & \iff & \Lambda_{\rho}\bigg(\Delta_{\omega}(\varphi_1-\varphi_2)\,\omega - \frac{1}{n}\,\Delta_{\rho}(\varphi_1-\varphi_2)\,\rho\bigg)=0\\
     & \iff & P_{\omega,\,\rho}(\varphi_1-\varphi_2)=0,\end{eqnarray}

\noindent where we have considered the operator

\begin{equation}\label{eqn:operator-diff-laplacians}P_{\omega,\,\rho}: = (\Lambda_{\rho}\omega)\,\Delta_{\omega} - \Delta_{\rho}.\end{equation}

\noindent Let $\rho_1,\dots, \rho_n>0$ be the eigenvalues of $\rho>0$ w.r.t. $\omega$. If we fix an arbitrary point $x\in X$ and choose local holomorphic coordinates $z_1, \dots , z_n$ about $x$ such that 

\vspace{1ex}

\hspace{8ex}$\omega(x)=\sum\limits_{j=1}^nidz_j\wedge d\bar{z}_j \hspace{2ex} \mbox{and} \hspace{2ex} \rho(x)=\sum\limits_{j=1}^n\rho_j(x)\,idz_j\wedge d\bar{z}_j,$

\vspace{1ex}

\noindent then $(\Lambda_{\rho}\omega)(x) = \sum\limits_{j=1}^n\frac{1}{\rho_j(x)}$ and $P_{\omega,\,\rho}(\varphi)(x) = \sum\limits_{j=1}^n\bigg(\sum\limits_{l\neq j}\frac{1}{\rho_l(x)}\bigg)\,\frac{\partial^2\varphi}{\partial z_j\partial\bar{z}_j}(x)$ for any real-valued $C^2$ function $\varphi$. This means that

\begin{equation}\label{eqn:P-formula}P_{\omega,\,\rho} = \Delta_{\tilde\lambda},\end{equation}

\noindent where $\tilde{\lambda}>0$ is the smooth $(1,\,1)$-form on $X$ whose eigenvalues w.r.t. $\omega$ are 

\begin{equation}\label{eqn:lambda-tilde-def}\frac{1}{\sum\limits_{l\neq j}\frac{1}{\rho_l}}>0, \hspace{1ex} j=1, \dots , n, \hspace{1ex} \mbox{hence}\hspace{1ex} \tilde{\lambda}(x):=\sum\limits_{j=1}^n\frac{1}{\sum\limits_{l\neq j}\frac{1}{\rho_l(x)}}\,idz_j\wedge d\bar{z}_j>0.\end{equation}

 We can actually give $\tilde{\lambda}$ an invariant expression. Let $\lambda$ be the smooth $(1,\,1)$-form intrinsically defined by

\begin{equation}\label{lambda_def}\lambda = \bigg[(\Lambda_{\rho}\omega)\,\frac{\omega^{n-1}}{(n-1)!} - \bigg(\frac{\omega^n}{\rho^n}\bigg)^{\frac{1}{n-1}}\,\star(\star\rho)^{\frac{1}{n-1}}\bigg]^{\frac{1}{n-1}}>0,\end{equation}

\noindent i.e. the $(n-1)^{st}$ root of a positive definite $(n-1,\,n-1)$-form. That this $(n-1,\,n-1)$-form is positive-definite follows from the calculation below showing it to be $\tilde{\lambda}^{n-1}$ multiplied by a positive function. We notice that using formulae (\ref{eqn:star_n-1-coord}) and (\ref{eqn:root-coeff}), we get

\begin{equation}\label{eqn:star-star-rho}\star(\star\rho)^{\frac{1}{n-1}} = \frac{1}{((n-1)!)^{\frac{1}{n-1}}}\,(\rho_1\dots\rho_n)^{\frac{1}{n-1}}\,\sum\limits_{j=1}^n\frac{1}{\rho_j}\,\widehat{idz_j\wedge d\bar{z}_j} \hspace{3ex} \mbox{at}\,\,x,\end{equation}

\begin{eqnarray}\label{eqn:lambda-coord}\nonumber\tilde{\lambda}^{n-1} & = & (n-1)!\,(\det\tilde{\lambda})\,\sum\limits_{j=1}^n\bigg(\sum\limits_{l\neq j}\frac{1}{\rho_l}\bigg)\,\widehat{idz_j\wedge d\bar{z}_j}\\
\nonumber & = & (n-1)!\,(\det\tilde{\lambda})\,\bigg(\sum\limits_{l=1}^n\frac{1}{\rho_l}\bigg)\,\sum\limits_{j=1}^n\widehat{idz_j\wedge d\bar{z}_j} - (n-1)!\,(\det\tilde{\lambda})\,\sum\limits_{j=1}^n\frac{1}{\rho_j}\,\widehat{idz_j\wedge d\bar{z}_j}\\
\nonumber & \stackrel{(a)}{=} & (n-1)!\,(\det\tilde{\lambda})\, \bigg[(\Lambda_{\rho}\omega)\,\frac{\omega^{n-1}}{(n-1)!} - (\rho_1\dots\rho_n)^{-\frac{1}{n-1}}\,\star(\star\rho)^{\frac{1}{n-1}}\bigg]\\
\nonumber & = & (n-1)!\,(\det\tilde{\lambda})\,\bigg[(\Lambda_{\rho}\omega)\,\frac{\omega^{n-1}}{(n-1)!} - \bigg(\frac{\rho^n}{\omega^n}\bigg)^{-\frac{1}{n-1}}\,\star(\star\rho)^{\frac{1}{n-1}}\bigg]\\
       & = & (n-1)!\,(\det\tilde{\lambda})\,\lambda^{n-1} \hspace{3ex} \mbox{at}\,\,x,\end{eqnarray}

\noindent where $(a)$ follows from (\ref{eqn:star-star-rho}). Taking determinants, we get $(\det\tilde{\lambda})^{n-1} = ((n-1)!)^n\,(\det\tilde{\lambda})^n\,(\det\lambda)^{n-1}$ at $x$, i.e. $\det\tilde{\lambda} = \frac{1}{((n-1)!)^n\,(\det\lambda)^{n-1}}$ at $x$. Thus (\ref{eqn:lambda-coord}) translates to $\tilde{\lambda} = (1/\det\lambda)\,\lambda$ at $x$. Now, in the chosen local coordinates, $\det\lambda = \lambda^n/\omega^n$ at $x$. Since $x\in X$ is arbitrary, the intrinsic shape of (\ref{eqn:lambda-coord}) is

\begin{equation}\label{eqn:lambda-tilde-lambda}\tilde{\lambda} = \frac{1}{(n-1)!}\,\frac{1}{\lambda^n/\omega^n}\,\lambda  \hspace{3ex} \mbox{on}\,\,X.\end{equation}

\noindent Combined with (\ref{lambda_def}), this is the sought-after invariant expression for $\tilde{\lambda}$. 

 Thanks to (\ref{eqn:P-formula}) and (\ref{eqn:lambda-tilde-lambda}), (\ref{gamma_1,2-eq2}) translates to the equivalences

$$\gamma_1^n = \gamma_2^n \iff \Delta_{\tilde\lambda}(\varphi_1 - \varphi_2)=0 \iff \Delta_{\lambda}(\varphi_1-\varphi_2)=0.$$

 By the maximum principle, the condition $\Delta_{\lambda}(\varphi_1-\varphi_2)=0$ on the compact manifold $X$ implies that $\varphi_1 - \varphi_2$ is constant on $X$. The proof of Proposition \ref{Prop:uniqueness-second} is complete.  \hfill $\Box$

\vspace{3ex}

 The uniqueness of solutions to equation $(\star)$ subject to (\ref{eqn:conditions}) is proved in the same way. The new terms are all of the first order, so they do not disturb in any way the ellipticity of the operators involved and the application of the maximum principle.

\begin{The}\label{The:uniqueness} Let $(X,\,\omega)$ be a compact Hermitian manifold, $\mbox{dim}_{\C}X=n\geq 2$. Suppose that for real-valued $C^{\infty}$ functions $\varphi_1$ and $\varphi_2$ on $X$ we have $\omega^{n-1} + i\partial\bar\partial\varphi_l\wedge\omega^{n-2} + \frac{i}{2}\,(\partial\varphi_l\wedge\bar\partial\omega^{n-2} - \bar\partial\varphi_l\wedge\partial\omega^{n-2}) > 0$ (for $l=1,2$) and

\begin{eqnarray}\nonumber\label{eqn:det-equality-general}\bigg[\bigg(\omega^{n-1} & + & i\partial\bar\partial\varphi_1\wedge\omega^{n-2} + \frac{i}{2}\,(\partial\varphi_1\wedge\bar\partial\omega^{n-2} - \bar\partial\varphi_1\wedge\partial\omega^{n-2})\bigg)^{\frac{1}{n-1}}\bigg]^n\\ 
\nonumber & = & \bigg[\bigg(\omega^{n-1} + i\partial\bar\partial\varphi_2\wedge\omega^{n-2} + \frac{i}{2}\,(\partial\varphi_2\wedge\bar\partial\omega^{n-2} - \bar\partial\varphi_2\wedge\partial\omega^{n-2})\bigg)^{\frac{1}{n-1}}\bigg]^n.\end{eqnarray}

\noindent Then the function $\varphi_1 - \varphi_2$ is constant on $X$.

\end{The}

\noindent {\it Proof.} For $l=1, 2$, we consider the positive definite $(1,\,1)$-forms

\vspace{1ex}

$\gamma_l:=\bigg(\omega^{n-1} + i\partial\bar\partial\varphi_l\wedge\omega^{n-2} + \frac{i}{2}\,(\partial\varphi_l\wedge\bar\partial\omega^{n-2} - \bar\partial\varphi_l\wedge\partial\omega^{n-2})\bigg)^{\frac{1}{n-1}}$

\vspace{1ex}

\noindent and the positive definite $(n-1,\, n-1)$-form

\vspace{1ex}

\hspace{15ex} $\Omega:=\sum\limits_{p=1}^n\bigg(\star(\gamma_1^{n-1})\bigg)^{n-p}\wedge\bigg(\star(\gamma_2^{n-1})\bigg)^{p-1}.$

\vspace{1ex}

\noindent If we set $\rho:=\Omega^{\frac{1}{n-1}}>0$, arguing as in the proof of Proposition \ref{Prop:uniqueness-second}, we find that the identity $\gamma_1^n=\gamma_2^n$ on $X$ is equivalent to

$$\bigg[\star\bigg(i\partial\bar\partial(\varphi_1-\varphi_2)\wedge\omega^{n-2}\bigg) + \frac{i}{2}\,\star\bigg(\partial(\varphi_1-\varphi_2)\wedge\bar\partial\omega^{n-2} - \bar\partial(\varphi_1-\varphi_2)\wedge\partial\omega^{n-2}\bigg)\bigg]\wedge\rho^{n-1}=0,$$ 

\noindent which, in turn, is found as in the proof of Proposition \ref{Prop:uniqueness-second} to be equivalent to

\begin{equation}\label{eqn:lambda-Q-tilde}(\Delta_{\tilde{\lambda}} + \tilde{Q})(\varphi_1-\varphi_2) = 0,\end{equation}

\noindent where $\tilde{\lambda}$ is the positive definite $(1,\,1)$-form defined by the eigenvalues $\rho_1,\dots,\rho_n$ $>0$ of $\rho$ w.r.t. $\omega$ through the same formula as in (\ref{eqn:lambda-tilde-def}) and $\tilde{Q}$ is the first-order operator defined on functions by

$$\tilde{Q}(\varphi):=\frac{i}{2}\,\star\bigg(\partial\varphi\wedge\bar\partial\omega^{n-2} - \bar\partial\varphi\wedge\partial\omega^{n-2}\bigg)\wedge\rho^{n-1}/\omega^n.$$

\noindent We can still define the $(1,\,1)$-form $\lambda>0$ intrinsically on $X$ by formula (\ref{lambda_def}) (only $\rho$ is different now) and it is still related to $\tilde{\lambda}$ by (\ref{eqn:lambda-tilde-lambda}). Setting 

\vspace{1ex}

\noindent $Q:=((n-1)!)^{\frac{1}{n-1}}\,\bigg(\prod\limits_{j=1}^n\bigg(\sum\limits_{l\neq j}\frac{1}{\rho_l}\bigg)\bigg)^{-\frac{1}{n-1}}\,\tilde{Q} = ((n-1)!)^{\frac{1}{n-1}}\,(\det\tilde\lambda)^{\frac{1}{n-1}}\,\tilde{Q} = \frac{1}{(n-1)!}\,\frac{1}{\det\lambda}\,\tilde{Q}$  \hspace{2ex} (see (\ref{eqn:lambda-coord})),

\vspace{1ex}

\noindent we see that (\ref{eqn:lambda-Q-tilde}) is equivalent to $(\Delta_{\lambda} + Q)(\varphi_1 - \varphi_2)=0$ on $X$. Since there are no zero-order terms in the second-order elliptic operator $\Delta_{\lambda} + Q$ and $X$ is compact, we conclude by the maximum principle that $\varphi_1 - \varphi_2$ is constant.  \hfill $\Box$

\section{The linearisation of equation $(\star)$}\label{section:linearisation}

 We will follow the analogy with the classical Calabi-Yau equation. We fix arbitrary $k\in\N$ ($k\geq 2$) and $0<\alpha<1$, and consider the open subset

$$U:=\{\varphi\in C^{k,\,\alpha}(X)\,\slash\,\omega^{n-1} + i\partial\bar\partial\varphi\wedge\omega^{n-2} + \frac{i}{2}\,(\partial\varphi\wedge\bar\partial\omega^{n-2} - \bar\partial\varphi\wedge\partial\omega^{n-2})>0\}$$

\noindent of the space $C^{k,\,\alpha}(X)$ of real functions on $X$ of H\"older class $C^{k,\,\alpha}$. We will calculate the differential at an arbitrary $\varphi\in U$ of the map $C:U\rightarrow C^{k-2,\,\alpha}(X)$,

$$C(\varphi) = \frac{\bigg[\bigg(\omega^{n-1} + i\partial\bar\partial\varphi\wedge\omega^{n-2} + \frac{i}{2}\,(\partial\varphi\wedge\bar\partial\omega^{n-2} - \bar\partial\varphi\wedge\partial\omega^{n-2})\bigg)^{\frac{1}{n-1}}\bigg]^n}{\omega^n}.$$

\noindent Let $\gamma>0$ be the smooth $(1,\,1)$-form on $X$ such that $\gamma^{n-1} = \omega^{n-1} + i\partial\bar\partial\varphi\wedge\omega^{n-2} + \frac{i}{2}\,(\partial\varphi\wedge\bar\partial\omega^{n-2} - \bar\partial\varphi\wedge\partial\omega^{n-2}):=\Lambda >0$. We will prove the following

\begin{Prop}\label{Prop:linearisation} For every $\varphi\in U$, the differential of $C$ at $\varphi$ calculated at an arbitrary $h\in C^{k,\,\alpha}(X)$ is given by the formula

\begin{equation}\label{eqn:linearisation}\nonumber C(\varphi)^{-1}(d_{\varphi}C)(h) =  \frac{1}{(n-1)^2}\,\bigg(\frac{\mbox{tr}_{\omega}\gamma}{\gamma^n/\omega^n}\,\Delta_{\omega}h - (n-1)!\,\Delta_{\star_{\omega}\Lambda}h\bigg) + \mbox{first order terms},\end{equation}

\noindent where the first order terms are $\frac{1}{n-1}\,\mbox{tr}_{\gamma^{n-1}}\bigg(\frac{i}{2}\,(\partial h\wedge\bar\partial\omega^{n-2} - \bar\partial h\wedge\partial\omega^{n-2})\bigg).$ 

\end{Prop}

 The rest of this section will be devoted to the proof of Proposition \ref{Prop:linearisation}. Using (\ref{eqn:determinants}), we see that $\log C(\varphi) = \frac{1}{n-1}\,\log \det(\gamma^{n-1}) + \frac{n}{n-1}\,\log(n-1)! - \log\det(\omega),$ where $\det(\gamma^{n-1})$ (resp. $\det(\omega)$) denotes the determinant of the coefficient matrix of $\gamma^{n-1}$ (resp. $\omega$) in local coordinates. Using the standard formula $(\log\det A)' = \mbox{tr}(A^{-1}A')$ and the fact that $\log\det(\omega)$ and $\frac{n}{n-1}\,\log(n-1)!$ do not depend on $\varphi$, we get\!\!:

\begin{eqnarray}\label{eqn:lin-eq1}\nonumber C(\varphi)^{-1}(d_{\varphi}C)(h) & = & \frac{1}{n-1}\,\mbox{tr}_{\gamma^{n-1}}\bigg(i\partial\bar\partial h\wedge\omega^{n-2}\bigg)\\
    & + & \frac{1}{n-1}\,\mbox{tr}_{\gamma^{n-1}}\bigg(\frac{i}{2}\,(\partial h\wedge\bar\partial\omega^{n-2} - \bar\partial h\wedge\partial\omega^{n-2})\bigg).\end{eqnarray}

\noindent We will now transform the first term in the r.h.s. above (i.e. the principal part of $C(\varphi)^{-1}d_{\varphi}C$). Setting as usual $\gamma_{n-1}:=\gamma^{n-1}/(n-1)!$, we get

\begin{eqnarray}\label{eqn:lin-eq2}\nonumber\frac{1}{n-1}\,\mbox{tr}_{\gamma^{n-1}}\bigg(i\partial\bar\partial h\wedge\omega^{n-2}\bigg) & = & \frac{(n-2)!}{(n-1)\,(n-1)!}\,\mbox{tr}_{\gamma_{n-1}}\bigg(i\partial\bar\partial h\wedge\frac{\omega^{n-2}}{(n-2)!}\bigg)\\
   & \stackrel{(a)}{=} & \frac{1}{(n-1)^2}\,\mbox{tr}_{\gamma}\bigg(\star_{\gamma}\bigg(i\partial\bar\partial h\wedge\frac{\omega^{n-2}}{(n-2)!}\bigg)\bigg)\end{eqnarray}

\noindent where identity $(a)$ has followed from (\ref{eqn:trace-star}) applied with $\gamma$ in place of $\omega$.

\begin{Lem}\label{Lem:trace-star-2metrics} For any Hermitian metrics $\omega,\gamma>0$ and any smooth real $(n-1,\,n-1)$-form $\Gamma$ on $X$, we have

\begin{equation}\label{eqn:trace-star-2metrics}\mbox{tr}_{\gamma}(\star_{\gamma}\Gamma) = \frac{1}{\gamma^n/\omega^n}\,\bigg\langle\gamma,\,\star_{\omega}\Gamma\bigg\rangle_{\omega}.\end{equation}

\end{Lem}

\noindent {\it Proof.} Fix an arbitrary point $x\in X$ and local holomorphic coordinates $z_1,\dots , z_n$ centred on $x$ such that

$$\omega(x)=\sum\limits_{j=1}^nidz_j\wedge d\bar{z}_j, \hspace{2ex} \gamma(x)=\sum\limits_{j=1}^n\gamma_j\,idz_j\wedge d\bar{z}_j, \hspace{2ex} \Gamma(x)=\sum\limits_{j,\,k=1}^n\Gamma_{j\bar{k}}\,\widehat{idz_j\wedge d\bar{z}_k},$$

\noindent with $\gamma_j>0$, $\Gamma_{j\bar{j}}\in\R$. Since $\langle dz_j\wedge d\bar{z}_j,\, dz_k\wedge d\bar{z}_k\rangle_{\gamma} = \frac{\delta_{jk}}{\gamma_j^2}$ at $x$, we get at $x$\!\!:

\vspace{1ex}

\noindent $\langle\widehat{dz_j\wedge d\bar{z}_j},\, \widehat{dz_k\wedge d\bar{z}_k}\rangle_{\gamma} = \delta_{jk}\,\frac{\gamma_j^2}{\gamma_1^2\cdots\gamma_n^2},$ so $\star_{\gamma}\widehat{(i dz_j\wedge d\bar{z}_j)} = \frac{\gamma_j^2}{\gamma_1\cdots\gamma_n}\,i dz_j\wedge d\bar{z}_j.$

\vspace{1ex}

\noindent It follows that $(\star_{\gamma}\Gamma)(x) = \sum\limits_{j=1}^n\frac{\Gamma_{j\bar{j}}\,\gamma_j^2}{\gamma_1\cdots\gamma_n}\,i dz_j\wedge d\bar{z}_j + \sum\limits_{j\neq k}\xi_{j\bar{k}}\,i dz_j\wedge d\bar{z}_k$, where the non-diagonal terms with coefficients denoted by $\xi_{j\bar{k}}$ can be computed explicitly but we will not do it here because they will disappear in the trace. 

 On the other hand, we have $\star_{\omega}\widehat{(i dz_j\wedge d\bar{z}_j)} = i dz_j\wedge d\bar{z}_j$ at $x$ by (\ref{eqn:star_n-1-coord}), hence $(\star_{\omega}\Gamma)(x) = \sum\limits_{j=1}^n\Gamma_{j\bar{j}}\,i dz_j\wedge d\bar{z}_j + \sum\limits_{j\neq k}\eta_{j\bar{k}}\,i dz_j\wedge d\bar{z}_k$. Therefore, we get

$$\mbox{tr}_{\gamma}(\star_{\gamma}\Gamma)(x) = \sum\limits_{j=1}^n\frac{\Gamma_{j\bar{j}}\,\gamma_j}{\gamma_1\cdots\gamma_n} = \frac{1}{\det\gamma(x)}\,\langle\gamma,\,\star_{\omega}\Gamma\rangle_{\omega}(x).$$

\noindent Since $\det\gamma(x)=\gamma^n/\omega^n(x)$ and $x$ is arbitrary, this proves the contention.  \hfill $\Box$

\vspace{2ex}

\noindent {\it End of proof of Proposition \ref{Prop:linearisation}.}

\vspace{1ex}

 We see that (\ref{eqn:star-iddbar_phi-n-2}) reads $\star_{\omega}(i\partial\bar\partial h\wedge\omega^{n-2}/(n-2)!) = -i\partial\bar\partial h + (\Delta_{\omega}h)\,\omega$, so applying (\ref{eqn:trace-star-2metrics}) with $\Gamma = i\partial\bar\partial h\wedge\omega^{n-2}/(n-2)!$, the term featuring on the right side of (\ref{eqn:lin-eq2}) becomes

\begin{eqnarray}\label{eqn:lin-eq3}\nonumber\mbox{tr}_{\gamma}\bigg(\star_{\gamma}\bigg(i\partial\bar\partial h\wedge\frac{\omega^{n-2}}{(n-2)!}\bigg)\bigg) & = & \frac{1}{\gamma^n/\omega^n}\,\bigg(-\langle\gamma,\,i\partial\bar\partial h\rangle_{\omega} + (\Delta_{\omega}h)\,\langle\gamma,\,\omega\rangle_{\omega}\bigg)\\
    & = & \frac{\mbox{tr}_{\omega}\gamma}{\gamma^n/\omega^n}\,\Delta_{\omega}h - \bigg\langle\frac{\gamma}{\gamma^n/\omega^n},\,i\partial\bar\partial h\bigg\rangle_{\omega}.\end{eqnarray}

 We will now transform the last term in (\ref{eqn:lin-eq3}). Recall that $\gamma>0$ has been defined as the $(n-1)^{st}$ root of $\Lambda:=\gamma^{n-1} = \omega^{n-1} + i\partial\bar\partial\varphi\wedge\omega^{n-2} + \frac{i}{2}\,(\partial\varphi\wedge\bar\partial\omega^{n-2} - \bar\partial\varphi\wedge\partial\omega^{n-2})>0$. If at a given point $x\in X$ the local coordinates are chosen as in the proof of Lemma \ref{Lem:trace-star-2metrics}, then thanks to (\ref{eqn:root-coeff}) we have

\begin{equation}\label{eqn:power-root-again}\Lambda(x) = \sum\limits_{j=1}^n\Lambda_j\,\widehat{idz_j\wedge d\bar{z}_j} \hspace{2ex}\mbox{with}\hspace{2ex} \gamma_j=(n-1)!\,\frac{\det\gamma}{\Lambda_j}.\end{equation}

\noindent Thus it follows from (\ref{eqn:star_n-1-coord}) that the $(1,\,1)$-form $\star_{\omega}\Lambda>0$ is given at $x$ by $(\star_{\omega}\Lambda)(x) = \sum\limits_{j=1}^n\Lambda_j\,idz_j\wedge d\bar{z}_j$. Putting the various bits together, we get at $x$\!\!:

\begin{eqnarray}\nonumber\langle\gamma,\,i\partial\bar\partial h\rangle_{\omega} =  \sum\limits_{j=1}^n\gamma_j\,\frac{\partial^2h}{\partial z_j\partial\bar{z}_j} & = & (n-1)!\,(\det\gamma)\,\sum\limits_{j=1}^n\frac{1}{\Lambda_j}\,\frac{\partial^2h}{\partial z_j\partial\bar{z}_j}\\
\nonumber & = & (n-1)!\,(\det\gamma)\,\mbox{tr}_{\star_{\omega}\Lambda}(i\partial\bar\partial h),\end{eqnarray}

\noindent which in invariant terms translates to

 \begin{equation}\label{eqn:lin-eq4}\bigg\langle\frac{\gamma}{\gamma^n/\omega^n},\,i\partial\bar\partial h\bigg\rangle_{\omega} = (n-1)!\,\Delta_{\star_{\omega}\Lambda}h.\end{equation}

\noindent Proposition \ref{Prop:linearisation} follows by putting together (\ref{eqn:lin-eq1}), (\ref{eqn:lin-eq2}), (\ref{eqn:lin-eq3}) and (\ref{eqn:lin-eq4}).    \hfill $\Box$

\vspace{2ex}

 Part $(b)$ of Theorem \ref{The:equation-summing} appears now as an immediate consequence of Proposition \ref{Prop:linearisation}.

\begin{Cor}\label{Cor:linearisation} Set $\rho:=\star_{\omega}\Lambda>0$, a smooth $(1,\,1)$-form. Then, for every $\varphi\in U$, the principal part of $C(\varphi)^{-1}d_{\varphi}C$ is the second-order elliptic operator

\begin{equation}\label{eqn:linearisation}\nonumber\frac{(n-1)!}{(n-1)^2}\,\bigg((\Lambda_{\rho}\omega)\,\Delta_{\omega} - \Delta_{\rho}\bigg) = \frac{(n-2)!}{n-1}\,P_{\omega,\,\rho} = \frac{(n-2)!}{n-1}\,\Delta_{\tilde\lambda},\end{equation}

\noindent with $P_{\omega,\,\rho}$ (resp. $\tilde\lambda$) defined in terms of $\omega$ and $\rho$ by formula (\ref{eqn:operator-diff-laplacians}) (resp. (\ref{eqn:lambda-tilde-def})).

\end{Cor}

\noindent {\it Proof.} Fix a point $x\in X$. We keep the notation and the choice of local coordinates $z_1,\dots , z_n$ centred on $x\in X$ of the proof of Proposition \ref{Prop:linearisation}. At $x$, we have

$$\frac{\mbox{tr}_{\omega}\gamma}{\gamma^n/\omega^n} = \frac{1}{\det\gamma}\,\sum\limits_{j=1}^n\gamma_j \stackrel{(a)}{=} (n-1)!\,\sum\limits_{j=1}^n\frac{1}{\Lambda_j} = (n-1)!\,\mbox{tr}_{\star_{\omega}\Lambda}\omega = (n-1)!\,\Lambda_{\rho}\omega.$$

\noindent where $(a)$ has followed from (\ref{eqn:power-root-again}). Combined with the conclusion of Proposition \ref{Prop:linearisation} and with (\ref{eqn:P-formula}), this proves the contention.   \hfill  $\Box$

\vspace{3ex}

\noindent {\bf References.} \\

\noindent [Dem92]\, J.-P. Demailly --- {\it Regularization of Closed Positive Currents and Intersection Theory} --- J. Alg. Geom., {\bf 1} (1992), 361-409.

\vspace{1ex}

\noindent [DP04]\, J.-P. Demailly, M. Paun --- {\it Numerical Charaterization of the K\"ahler Cone of a Compact K\"ahler Manifold} --- Ann. of Math. (2) {\bf 159(3)} (2004) 1247-1274.

\vspace{1ex}

\noindent [FWW10a]\, J. Fu, Z. Wang, D. Wu --- {\it Form-type Calabi-Yau Equations} --- Math. Res. Lett. {\bf 17} (2010), no. 5, 887-903. 

\vspace{1ex}

\noindent [FWW10b]\, J. Fu, Z. Wang, D. Wu --- {\it Form-type Calabi-Yau Equations on K\"ahler Manifolds of Nonnegative Orthogonal Bisectional Curvature} -- arXiv e-print math.DG/1010.2022.

\vspace{1ex}

\noindent [Fri91]\, R. Friedman --- {\it On Threefolds with Trivial Canonical Bundle} --- Complex Geometry and Lie Theory (Sundance, UT, 1989) 103-134, Proc. Sympos. Pure Math., {\bf 53} , Amer. Math. Soc, Providence R.I. 1991.

\vspace{1ex}

\noindent [FLY12]\, J.Fu, J.Li, S.-T. Yau --- {\it Balanced Metrics on Non-K\"ahler Calabi-Yau Threefolds} --- J. Diff. Geom. {\bf 90} (2012) 81-129.

\vspace{1ex}

\noindent [Gau77]\, P. Gauduchon --- {\it Le th\'eor\`eme de l'excentricit\'e nulle} --- C.R. Acad. Sc. Paris, S\'erie A, t. {\bf 285} (1977), 387-390.

\vspace{1ex}

\noindent [IP12]\, S. Ivanov, G. Papadopoulos --- {\it Vanishing Theorems on $(l/k)$-strong K\"ahler Manifolds with Torsion} --- arXiv e-print DG 1202.6470v1.

\vspace{1ex}

\noindent [KS60]\, K. Kodaira, D.C. Spencer --- {\it On Deformations of Complex Analytic Structures, III. Stability Theorems for Complex Structures} --- Ann. of Math. {\bf 71}, no.1 (1960), 43-76

\vspace{1ex}

\noindent [LT96]\, P. Lu, G. Tian --- {\it Complex Structures on Connected Sums of $S^3\times S^3$} --- Manifolds and Geometry (Pisa, 1993), 284-293, Sympos. Math., XXXVI, Cambridge Univ. Press, Cambridge, 1996.

\vspace{1ex}

\noindent [Mic82]\, M. L. Michelsohn --- {\it On the Existence of Special Metrics in Complex Geometry} --- Acta Math. {\bf 149} (1982), no. 3-4, 261-295.

\vspace{1ex}

\noindent [Pop09]\,D. Popovici --- {\it Deformation Limits of Projective Manifolds\!\!: Hodge Numbers and Strongly Gauduchon Metrics} --- Invent. Math. {\bf 194} (2013), 515-534.

\vspace{1ex}

\noindent [Pop13a]\, \, D. Popovici ---{\it Holomorphic Deformations of Balanced Calabi-Yau $\partial\bar\partial$-Manifolds} --- arXiv e-print math.AG/1304.0331v1

\vspace{1ex}

\noindent [Sch07]\, M. Schweitzer --- {\it Autour de la cohomologie de Bott-Chern} --- arXiv e-print math.AG/0709.3528v1.

\vspace{1ex}

\noindent [ST10]\, J. Streets, G. Tian --- {\it A Parabolic Flow of Pluriclosed Metrics} --- Int. Math. Res. Notices, {\bf 16} (2010) 3101-3133.

\vspace{1ex}

\noindent [Sul76]\, D. Sullivan --- {\it Cycles for the Dynamical Study of Foliated Manifolds and Complex Manifolds} --- Invent. Math. {\bf 36} (1976) 225 - 255.

\vspace{1ex}

\noindent [TW13a]\, V. Tosatti, B. Weinkove --- {\it The Monge-Ampere Equation for $(n-1)$-plurisubharmonic Functions on a Compact K\"ahler Manifold} --- arXiv e-print math.DG/1305.7511. 

\vspace{1ex}

\noindent [TW13b]\, V. Tosatti, B. Weinkove --- {\it Hermitian Metrics, $(n-1, n-1)$-forms and Monge-Amp\`ere Equations} --- arXiv e-print math.DG/1310.6326.

\vspace{1ex}

\noindent [Voi02]\, C. Voisin --- {\it Hodge Theory and Complex Algebraic Geometry. I.} --- Cambridge Studies in Advanced Mathematics, 76, Cambridge University Press, Cambridge, 2002.

\vspace{1ex}

\noindent [Wu06]\, C.-C. Wu --- {\it On the Geometry of Superstrings with Torsion} --- thesis, Department of Mathematics, Harvard University, Cambridge MA 02138, (April 2006).

\vspace{3ex}

\noindent Institut de Math\'ematiques de Toulouse, Universit\'e Paul Sabatier,

\noindent 118 route de Narbonne, 31062 Toulouse, France

\noindent Email\!: popovici@math.univ-toulouse.fr

\end{document}